\documentclass[a4paper,10pt]{amsart}
\usepackage{mathrsfs}
\usepackage[all]{xy}
\usepackage{calrsfs}
\usepackage{times}
\usepackage[scaled=0.92]{helvet}
\usepackage[greek,english]{babel}
\usepackage[iso-8859-7]{inputenc}

\usepackage{amsfonts}
\usepackage{amssymb}
\usepackage[usenames]{color}

\usepackage[usenames]{color}
\newtheorem{theorem}{Theorem}   

\newtheorem{lemma}[theorem]{Lemma}
\newtheorem{proposition}[theorem]{Proposition}
\newtheorem{corollary}[theorem]{Corollary}

\theoremstyle{definition}
\newtheorem{definition}[theorem]{Definition}
\newtheorem{remark}[theorem]{Remark}
\newtheorem{example}[theorem]{Example}

\definecolor{MyDarkBlue}{rgb}{0,0.08,0.60}

\newcommand{\A}{{ \rm Aut }}

\newcommand{\F}{{\mathbb{F}}}

\newcommand{\Z}{{\mathbb{Z}}}

\newcommand{\lf}{\left\lfloor}
\newcommand{\rf}{\right\rfloor}

\renewcommand{\mod}{{\;\rm mod}}

\title{Automorphisms of Curves and Weierstrass semigroups}
\author{Sotiris Karanikolopoulos \and Aristides Kontogeorgis}
\address{
Department of Mathematics, University of the \AE gean, 83200 Karlovassi, Samos,
Greece
}
\email{kontogar@aegean.gr,mathm03005@aegean.gr}

\date{\today}

\begin{document}
\bibliographystyle{amsplain}

\begin{abstract}
The relation of the Weierstrass semigroup with several 
invariants of a curve is studied.  
For Galois covers of curves with group $G$ we introduce a
new filtration of the group decomposition subgroup of $G$.
The relation to the ramification filtration is 
investigated in the case of cyclic covers.
We relate our results to invariants defined by Boseck
and we study the one point ramification case. 
We also give applications to Hasse-Witt invariant and 
symmetric semigroups.

\end{abstract} 

\thanks{{\bf keywords:} Automorphisms, Curves, Differentials, Numerical Semigroups,
Hasse-Witt invariant, maximal curves. 
 {\bf AMS subject classification} 14H37,11G20}

\maketitle 
\section{Introduction}
Let $X$ be a projective nonsingular curve of genus $g\geq 2$ defined over an algebraic closed field $k$
of characteristic $p>0$. Let $G$ be a subgroup of the automorphism 
group $\A(X)$ of $X$ and let $G(P)$ be the subgroup of automorphisms stabilizing  a point $P$ on $X$. 
Consider also the Weierstrass semigroup attached at  the point $P$.
In characteristic zero there are results \cite{MorrisonPinkham} relating the structure of the Weierstrass 
semigroup at $P$ 
to the subgroup $G(P)$. In this paper we will try to prove analogous results in the positive characteristic case.
It turns out that there is a close connection of the automorphism group, the Weierstrass semigroup and 
other invariants 
of the curve like the rank of the Hasse-Witt matrix. Aim of this article is to show how various invariants 
of a curve are  encoded in the Weierstrass semigroup. 

Section \ref{section-2} is devoted to connections of the theory of decomposition groups $G(P)$ 
to the theory of Weierstrass semigroups.  
Both the theory of the Weierstrass semigroups and the decomposition groups $G(P)$ are more difficult 
when $p>0$. 
In characteristic zero it is known that $G(P)$ is always a cyclic group, while when $p>0$ and
 $p$ divides $|G(P)|$ the 
group $G(P)$ is no more cyclic and admits the 
following ramification filtration:
\begin{equation*} 
 G(P)=G_0(P)   \supset G_1(P) \supset G_2(P) \supset \ldots,
\end{equation*}
Recall that the groups  $G_i(P)$ are defined as 
  $G_i(P)=\{\sigma \in G(P): v_P(\sigma(t)-t )\geq i+1 \}$, for a local uniformizer 
 $t$  at $P$ and  $v_P$ is the corresponding valuation. Notice that $G_1(P)$ is the $p$-part of $G(P)$. 
For every point $P$  of the curve $X$  
we consider the sequence of $k$-vector 
spaces 
\begin{equation*} 
k=L(0)=L(P)= \cdots  =L((i-1)P) < L( iP) \leq  \cdots \leq  L( (2g-1)P ),
\end{equation*}
where
\[
L(iP):=\{ f \in k(X)^*: {\rm div}(f) +iP \geq 0 \} \cup \{0\}.
\]
We will write $\ell( D)= \dim_k L(D)$.
An integer $i$ will be called  a pole number if there is a function 
$f \in k(X)^*$ so that $\mathrm{div}_\infty(f)=iP$ or 
equivalently $\ell\big( (i-1)P \big)+1=\ell\big(i P\big)$. 
The set of pole numbers at $P$  form a semigroup $H(P)$ which is called the 
Weierstrass semigroup at $P$. 
It is known that there are exactly $g$ pole numbers that are smaller or equal to $2g-1$ and that every 
integer $i \geq 2g-1$ is in the Weierstrass semigroup, see \cite[I.6.7]{StiBo}.  

\begin{remark}In  characteristic zero for all but finite points (the so called Weierstass points),
 the generic situation 
for the gaps at $P$ is  the set $\{1,\ldots,g\}$. 
This is not correct in positive characteristic,
where the generic set of gaps might be different than the set   $\{1,\ldots,g\}$.
\end{remark}

\begin{remark}
 The Weierstrass semigroup seems to be defined locally at a given point $P$. 
But the condition on the existence of function with only pole at $P$ is a global condition and
lot of global invariants like the genus or the Hasse-Witt matrix \cite{Stoehr-Vianna} are encoded in it.  
\end{remark}

In section \ref{sectionRamGroups} we recall the relation of the representation of the
$p$-part of the decomposition group into the Riemann-Roch spaces introduced in \cite{KontoZ}.
In section \ref{Jumps} we introduce the {\em representation filtration} that shares many properties with 
the ramification filtration and affects the structure of the Weierstrass semigroup. 
In propositions \ref{10prop} and \ref{11prop} we characterize the jumps of the representation 
filtration  in terms of the Weierstrass semigroup and in terms of the intermediate field extensions.
 Valentini
and Madan  \cite[lemma 2]{Valentini-MadanZ} did a  similar computation
for the gaps
of the ramification filtration for cyclic groups using the theory of  Witt
vectors, see \cite{vm2},\cite{Schmid36}.

Section \ref{cyc-section} is devoted to the  description of the Weierstrass semigroup when 
$G_1(P)$ is a cyclic group. This is an important case to understand,  because of its 
simplicity and because of  Oort conjecture regarding the deformation of a curve with 
automorphisms from characteristic $p$ to characteristic $0$. 
These results should be seen as the characteristic $p$ analogon of the work of 
Morrison-Pinkham \cite{MorrisonPinkham} on Weierstrass semigroups for 
Galois Weierstrass points.

In section \ref{hol-diff} we still study the semigroups that can appear in the cyclic group 
case but using the Boseck invariants. These are invariants that were first introduced by Boseck 
in \cite{boseck} in the study of basis of holomorphic differentials of curves that 
are cyclic extensions of the rational function fields. 
Boseck in his very important and interesting article studied the Weierstrass semigroups in this case.
Unfortunately, he didn't notice that there is one case that needs more attention, 
namely the case of a cover where only one point is ramified and the prime number 
is ``small''. What we mean by small will be clear in the sequel. For example for the 
case of Artin-Schreier curves small means that the characteristic is smaller than 
the conductor. The same problem is now in several places in the literature like 
\cite[p.235 remark (i)]{garciaa-s}, \cite[p. 170 in first 3 lines]{Val-Mad:80}, 
\cite[p. 3028, remark (i)]{garciaelab}. For the special case of Artin-Schreier 
extension the correct treatment is the one in lemma \ref{smallGaps}.

The curves that appear as covers of the projective line 
with only one ramified point are very interesting. This is because curves with 
large automorphisms groups appear this way \cite{StiI} and because of the 
Katz-Gabber compactification of local actions \cite[Thm 1.4.1]{Katz86},\cite[Cor. 1.9]{Gille00}. 

In section \ref{sec-semigroups} we study the maximal gap in the semigroup and we  give 
bounds for it, in terms of several invariants of the curve. When this maximal gap equals to
$2g-1$ at some point $P$, then the Weierstrass semigroup is called \textit{symmetric} 
and that point is a Weierstrass point, proposition \ref{genera}.
 The existence of a single point in the curve that 
has a symmetric Weierstrass semigroup determines the genus of the curve, proposition 
\ref{genera}. Many well known examples of curves  enjoy that property. 
In particular, we determine the cyclic totally ramified coverings of the rational function
 field with this symmetric  Weierstrass semigroups in lemmata \ref{cyc-symmetric} and \ref{cyc-symms}

Section \ref{sec-HasseWitt} is devoted to the dependence of the Cartier operator 
and the Hasse-Witt invariant to the theory of Weierstrass semigroups. 
We give some results about \textit{big actions}, lemma \ref{big11},
 corollary \ref{baoderG_2} and characterize some 
 non classical curves with respect to the canonical linear series and 
nilpotent Cartier operator, corollary \ref{nilpotentnonclas}. 
We also characterize all the maximal curves over $\mathbb{F}_{q^2}$ 
with two generators for the Weierstrass semigroup at a $\mathbb{F}_{q^2}$--rational 
point, as curves with symmetric Weierstrass semigroups at this point, remark \ref{lastremark}.
 Finally we study further the connection between symmetric Weierstrass 
semigroups and maximal curves.

\section{Decompostion Groups $G_i(P)$ }
\label{section-2}
\subsection{Ramification Groups}
 \label{sectionRamGroups}

\begin{proposition} \label{fai-rep-lemma}
If $g\geq 2$ and $p\neq 2,3$ then there is at least 
one pole number $m \leq 2g-1$ not divisible by the characteristic $p$.
        Let $1<m\leq 2g-1$ be the smallest  pole number not divisible by the characteristic. 
There is  a faithful representation 
\begin{equation} \label{fai-rep}
\rho: G_1(P) \rightarrow \mathrm{GL} \big( L(mP) \big)
\end{equation}
\end{proposition}
\begin{proof}
\cite[lemmata  2.1,2.2]{KontoZ}
\end{proof}

\begin{proposition}
A basis for the vector space $L(mP)$ is given by 
\[
\left\{1,\frac{u_i}{t^{m_i}},\frac{1}{t^m}:\mbox{ where } 1<i<r, p \mid m_i \mbox{ and } u_i 
\mbox{ are certain  units} \right\}.
\]
With respect to  this basis, an element  $\sigma \in G_1(P)$ acts on $1/t^m$ by 
\[
\sigma\frac{1}{t^m}=\frac{1}{t^m}+ \sum _{i=1}^r c_i(\sigma) 
 \frac{u_i}{t^{m_i}}.
\]
The action on the local uniformizer is given by 
\begin{equation} 
\sigma(t)=t-\frac{1}{m} \sum_{i=1}^r c_i(\sigma) u_i 
t^{m-m_i+1}+\cdots ,
\end{equation}
\end{proposition}
and this implies 
\begin{proposition} \label{gap-comp}
Let $P$ be a wild ramified  point on the curve $X$ and let 
 \[\rho: G_1(P) \rightarrow \mathrm{GL}_{\ell(mP)} (k)\] be the corresponding  faithful 
representation we considered  in  Lemma \ref{fai-rep-lemma}. 
Let $m=m_r>m_{r-1} > \cdots > m_0=0$ be the pole numbers at $P$ that are 
$ \leq m$. If $G_{i}(P)  >  G_{i+1}(P)$  then   $ i=m-m_k$, 
for some pole number $m_k$. 
\end{proposition}
\begin{remark}
Notice that in this article we enumerate the pole numbers in terms of an increasing function 
$i\mapsto m_i$. In \cite{KontoZ}  the enumeration is in terms of a decreasing function. 
\end{remark}
\subsection{Jumps in the ramification filtration and divisibility of the Weierstrass semigroup}
\label{Jumps}
Consider a Galois cover $\pi:X\rightarrow Y=X/G$ of algebraic curves,
and let $P$ be a fully ramified point of $X$. How are the Weierstrass semigroup 
sequences of $P$, and $\pi(P)$ related?
\begin{lemma} \label{up-down}
Let $k(X)$, $k(Y)=k(X)^G$  denote the function fields of the curves $X$ and 
$Y$ respectively. The morphisms 
\[
 N_G:k(X) \rightarrow k(Y) \mbox{ and } \pi^*:k(Y) \rightarrow k(X),
\]
 sending $f\in k(X)$ to $N_G(f)=\prod_{\sigma_\in G} f$ and $g\in k(Y)$ 
to $\pi^* g\in k(X)$ respectively, induce 
injections 
\[
N_G: H(P) \rightarrow H(Q)  \mbox{ and } \pi^*:H(Q) \stackrel{\times |G|}\longrightarrow H(P),
\]
 where $Q:=\pi (P)$.
\end{lemma}
\begin{proof}
For every element $f \in k(X)$ such that $(f)_\infty=mP$, 
the element $N_G(f)$ is a $G$-invariant element,  so it is in $k(Y)$. 
Moreover, the pole order of $N_G(f)$ seen as a function on $k(X)$ is 
$|G|\cdot m$. But since $P$ is fully ramified the valuation of $N_G(f)$ 
expressed in terms of the local uniformizer at $\pi(P)$ is just $-m$. 

On the other hand side an element $g\in k(Y)$ seen as an element 
of $k(X)$ by considering the pullback $\pi^*(g)$ has for the same reason 
valuation at $P$ multiplied by the order of $G$.
\end{proof}

\begin{remark}
 The condition of fully ramification is necessary in the above lemma. 
Indeed, if a point $Q \in Y$ has more than one elements in $\pi^{-1}(Q)$
then the pullback of $g$, such that $(g)_\infty=m Q$,  is supported 
on $\pi^{-1}(Q)$ and gives no information for the Weierstrass semigroup of 
any of the points $P \in \pi^{-1}(Q)$. 
\end{remark}

We will prove the following:
\begin{lemma} \label{div-inv}
 If an element $f$ such that $(f)_\infty = a P$ is invariant under the 
action of a subgroup $H < G_1(P)$, then $|H|$ divides $a$.
\end{lemma}
\begin{proof}
Write $f=u/t^a$ in terms of a local uniformizer $t$ at $P$ and a unit $u$. 
Since $f$ is invariant it is the pullback of a function $g \in k(X/H)$. 
If $t'$ is a local uniformizer at $Q=\pi(P)$ then $g$ is 
expressed as $g=u'/t'^b$. Since $t'=t^{|H|} v$ for some unit $v$, 
\cite[IV 2.2c]{Hartshorne:77} the desired result follows.
\end{proof}

\begin{definition}\label{repfil}
 Let  $0=m_0<m_1 < m_2 < \ldots < m_r$ be the sequence of pole numbers up to $m_r$. 
For each $ 0 \leq i \leq r$ we consider the representations
\[
 \rho_i :G_1(P) \rightarrow \mathrm{GL}(L(m_i P)).
\]
We form the decreasing sequence of groups:
\begin{equation} \label{rep-filt}
 G_1(P)=\mathrm{ker} \rho_0 \supseteq \mathrm{ker} \rho_1 \supseteq \mathrm{ker} \rho_2 
\supseteq \cdots \supseteq 
\mathrm{ker} \rho_r=\{1\}.
\end{equation}
We will cal this sequence of groups ``the representation'' filtration.
\end{definition}

Let $\sigma \in \mathrm{ker} \rho_i$. Then $\rho_{i+1}(\sigma)$ 
has the following form
\[
 \rho_{i+1}(\sigma)=
\begin{pmatrix}
 1 & 0 &  0 &\cdots   & 0 \\
 0 & 1 &   0     &\cdots  & 0  \\
 \vdots & 0 &\ddots &  \ddots&\vdots \\
0 & \cdots & 0 & 1  & 0 \\
a_{i+1,1}(\sigma) & a_{i+1,2}(\sigma) & \cdots & a_{i+1,i}(\sigma) & 1
\end{pmatrix}.
\]
Observe also that all functions 
$a_{i+1,\nu}:\mathrm{ker} \rho_i \rightarrow k$
are group homomorphisms into the additive group of the field $k$.
Notice that 
\[
 \mathrm{ker} \rho_{i+1}=\mathrm{ker} \rho_i \cap  \bigcap_{\nu=i+1}^i \mathrm{ker} a_{i+1,\nu}.
\]

\begin{lemma} \label{gallin}
The linear series  $|m_1 P|=\mathbb{P}(L(m_1 P))$ defines a map $X \rightarrow \mathbb{P}^1$ of degree $m_1$. 
This gives rise to an algebraic extension $F/k(f_1)$ of degree $m_1$, where $(f_1)_\infty=m_1P$. Let 
$F$ denote the function field of the curve $X$.
The extension $F/k(f_1)$ is Galois if and only if   $m_1=|\ker \rho_1|$, and in this 
case the Galois group is $\ker \rho_1$.
\end{lemma}
\begin{proof}
 Observe that $F \supseteq F^{\ker \rho_1} \supseteq k(f_1)$.
Notice that the degree of the extension $F/k(f_1)$ is $m_1$ therefore if 
$|\ker \rho_1|=m_1$ then $k(f_1)=F^{\ker \rho_1}$.

Assume now that the extension $F/k(f_1)$ is Galois with Galois group $A$. Since 
$F^{\ker \rho_1} \supseteq k(f_1)=F^A$ we have that $A \supseteq \ker \rho_1$. 
By definition $A=\{\sigma \in G_1(P): \sigma(f_1)=f_1\} \subseteq \ker \rho_1$.
\end{proof}

\begin{remark}
In characteristic zero  a point $P$ such that the covering map $X\rightarrow \mathbb{P}^1$ introduced by 
the linear system $|m_1 P|$ is Galois,
is called {\em Galois Weierstrass point}. 
The Weierstrass semigroup 
at a Galois Weierstrass point in characteristic zero is now well understood, \cite{MorrisonPinkham}.
Notice that in positive characteristic the condition $m_1\leq g$ does not imply 
that $P$ is a Weierstrass point, since the canonical linear system might be non 
classical. In section \ref{cyc-section}  we will restrict ourselves to cyclic totally ramified
 extensions of order 
$p^n$, $n\geq 2$ and then $P$ is a Weierstrass point 
 according to theorem 2 in \cite{vm2}. 
\end{remark}

Lemma \ref{div-inv} can be generalized as follows:
\begin{proposition}
The order $|\mathrm{ker} \rho_i |$ of the group $\mathrm{ker} \rho_i $
 divides $m_\nu$, for all $\nu \leq i$.
\end{proposition}
\begin{proof}
This comes out from lemma \ref{div-inv} since the elements $f_i$ such that
$(f_i)_\infty=m_i P$ are by definition $\ker \rho_i$-invariant.
\end{proof}

We recall the following definition:
\begin{definition}
Given a finite group $G$ the Frattini subgroup $\Phi(G)$ is defined to be the intersection of all
proper maximal subgroups of $G$. If $G$ is a $p$-group then 
the Frattini group is $\Phi(G)=G^p [G,G]$, and is the minimal normal subgroup $N$ of $G$
such that $G/N$ is elementary abelian. 
\end{definition}

Consider the first pole number $m_1$ in the Weierstrass semigroup and a function $f_1$ such that 
$(f_1)_\infty=m_1 P$. The image of the representation 
\[
 \rho:G_1(P) \rightarrow \mathrm{GL}(L(m_1P))
\]
is an elementary abelian group since the representation is $2$-dimensional, therefore the commutator and the Frattini subgroup 
of $G_1(P)$ are in the kernel of $\rho$. This combined with lemma \ref{div-inv} gives 
 us the following 
\begin{lemma} \label{div-Frat}
The order  $|\Phi(G)|$ of $\Phi(G)$ divides the first
pole number.
\end{lemma}

%
%
%
%
%
%
\begin{example}
In the work of C. Lehr, M. Matignon \cite{CL-MM}  a notion of {\em big action} 
is defined. Curves having a big action were studied further by M. Rocher and M. Matignon
\cite{Ma-Ro},\cite{Ro}.
A curve $X$ together with a subgroup $G$ of  the automorphism group of $X$ 
is called a big-action if $G$ is a $p$-group and 
\[
 \frac{|G|}{g} > \frac{2p }{p-1}.
\]
All big actions have the following property \cite{CL-MM}:
\end{example}
\begin{proposition}
Assume that $(X,G)$ is a big action.
 There is a unique  point $P$ of $X$ such that $G_1(P)=G$, the group $G_2(P)$ is not trivial and strictly contained in $G_1$ and the quotient $X/G_2(P)\cong \mathbb{P}^1$. Moreover, the 
group $G$ is an extension of groups
\[
 0 \rightarrow G_2(P) \rightarrow G=G_1(P) \stackrel{\pi}\longrightarrow (\Z/p\Z)^v \rightarrow 0.
\]
\end{proposition}
Since the curve $X/G_2(P)$ is of genus zero, the Weierstrass semigroup equals the 
semigroup of natural numbers, therefore $|G_2(P)|\cdot  \mathbb{Z}_+$ is a subsemigroup 
of the Weierstrass semigroup at the unique fixed point $P$.
\begin{lemma} \label{big11}
Let $(X,G)$ be a big action.
 The smallest non trivial pole number $m_1$  in $H(P)$ is $|G_2(P)|$.
\end{lemma}
\begin{proof}
Let $f_1$ be the function in the function field $k(X)$ such that
$(f_1)_\infty=m_1 P$. By \cite{Ma-Ro} $\mathrm{ker}\rho_1=G_2(P)$ and the result follows 
by lemma \ref{gallin}.
\end{proof}

\begin{theorem} \label{Lewittes}
 Let $F_0$ be an algebraic function field of genus $g$, defined over an algebraically closed field 
$k$ of characteristic $p> 0$. Select a place $P_0$ of $F_0$  and consider an Artin-Schreier extension $F/F_0$ 
given by $F=F_0(y)$, $y^{p^h}-y=G$, where $v_{P_0}(G)=-m<0$, $(m,p)=1$. We also assume 
that $G$ has no other poles at other places of $F_0$. 
Let $P$ be the unique place of $F$ that lies above $P_0$.
Assume that the set of gaps of $F_0$ at $P_0$ are given by 
\[
 \mathcal{G}(P_0)=
\left\{ \begin{array}{ll}
\{1=h_1<h_2<\cdots<h_{g_{F_0}}\leq 2g_{F_0}-1\}, & \mbox{ if } g_{F_0}\geq 1\\
\emptyset, & \mbox{ if } g_{F_0}=0 
\end{array}
\right. .
\]
Then the gaps $\mathcal{G}(P)$ at $P$ are given by 
\[
 \mathcal{G}(P)=\bigcup_{i=1}^{p^h-1} A_i,
\]
with 
\[
 A_i=\left\{mi-j p^h: 1\leq j \leq \lf \frac{mi}{p^h}\rf\right\} \bigcup \left\{ mi+p^hh_j:j=1,2,\dots,g_{F_0} \right\}.
 \]
\end{theorem}
\begin{proof}
 This is combination of theorems  1 and 5 of Lewittes, in \cite{Lewittes:89}.
Keep in mind that we are working over an algebraically closed field, 
therefore all places have degree 1. Also Lewittes proves his theorem for the more 
difficult case $k=\mathbb{F}_{q}$ but the proof over an algebraically closed 
field is the same and even simpler. 
\end{proof}
\begin{remark} \label{corM}
 The numbers $m,p^h$ are always pole numbers, since it is not possible to express
them as a gap given by theorem \ref{Lewittes}. This means that $m \Z_+ + p^h\Z_+$ 
is always included  in  $H(P)$. If the genus of $F_0$ is not $0$, then the inclusion 
$m\Z_++p^h \Z_+ \subseteq H(P)$ can be a strict inclusion, see remark \ref{two-AS}.
\end{remark}

\begin{proposition}\label{10prop}
 If $i$ is a  jump for the  representation filtration, i.e., 
 $\ker \rho_i \varsupsetneq \ker \rho_{i+1}$,
then  $m_{i+1} \not \in \langle m_1,\ldots, m_i\rangle_{\mathbb Z_+}$.
Also the number of jumps in the representation filtration equals the minimal number of generators of 
$H(P)$ up to $m_r$. 
\end{proposition}
\begin{proof}
 Let $f_i$ be elements in the function field of the curve such that 
$\mathrm{div}_\infty(f_i)=m_i P$. Let $f_1$ be the first such function 
with a unique pole at $P$ of order $m_1$. 

Suppose first  that 
$m_2=2m_1,m_3=3m_1,\ldots, m_\nu= \nu m_1$ are the first 
$\nu$ pole numbers.  Then a basis for the space $L(m_\nu P)$ 
consists of $\{f_1,f_1^2,\ldots,f_1^\nu\}$ and it is clear that 
 $\ker \rho_1=\ker \rho_2 =\cdots =\ker \rho_\nu$.
 
Now we treat the general case. If $m_{i+1}$ is in the semigroup $\langle m_1,\ldots, m_i\rangle_{\mathbb Z_+}$
generated by all $m_1,\ldots,m_i$ then $m_{i+1}=\sum_{\nu_i \in \mathbb{Z}_+} \nu_i m_i$
and 
\[
 f_{i+1}=\prod_{\nu_i \in \mathbb{Z}_+} f_i^{\nu_i}.
\]
 This implies that $\ker \rho_{i+1}=\ker \rho_i$.
\end{proof}

The next proposition characterize exactly the jumps of the representation filtration:
\begin{proposition} \label{11prop}
Denote by $X_i$ the curve $X/\ker \rho_i$ and by $F_i$ the corresponding function field.
  Supose that $X_1=\mathbb{P}^1$ and $F_1=k(f_1)$. 
The jumps at the representation 
filtration appear exactly at $i$ such that $F_i \neq F_{i+1}$. Moreover if 
{ $i_1,i_2,\ldots,i_n$
are  the }jumps of the representation filtration then 
 for all  $\ell$, with $1\leq \ell \leq n$ we have 
{$F_{i_\ell}=k(f_0,f_{i_1},\ldots, f_{i_\ell})$,} i.e. the jumps appear exactly at the integers $i$ such that 
the  function $f_i$ corresponding to the pole number $m_i$ gives rise to an element that can not 
be expressed as a rational function of $f_\nu$, $\nu<i$. 
\end{proposition}
\begin{proof}
We have the following picture
{
\[
 G_1(P)=\mathrm{ker} \rho_0 \supseteq
 \mathrm{ker} \rho_1  \supseteq \ldots \mathrm{ker} \rho_i\ldots  \supseteq \mathrm{ker}\rho_{r}=\{1\} 
\]
for the representation filtration and the strict inclusion holds for a subscript $i$ only when $i\in \{i_1,\ldots i_n\}$. In general $r+1\geq n$ with equality holds if all pole numbers up to $m_r$ are generators, or equivalently when all the inclusions above are strict.} We also have
\[
\xymatrix{
 X/G_1(P)\ar@{<->}[d] \ar[r] &  X/\mathrm{ker} \rho_1\ar@{<->}[d] \ar[r] & \ldots X/\mathrm{ker} \rho_i\ar@{<->}[d]\ldots  \ar[r]& X\ar@{<->}[d]\\
 F_0:=F_{r}^{G_1(P)} \ar@{-}[d]\subseteq &F_1:=k(f_1) \ar@{-}[d]\subseteq& F_i\ar@{-}[d]   &\subseteq  F_r \ar@{-}[d]\\
Q_0  \ar@{-}[r]  &Q_1 \ar@{.}[r]&   Q_i \ar@{. }[r]&  P
}
\]
where the first line corresponds to coverings, the second to function fields and the last to their totally ramified places.


Notice that $\ker \rho_r=\ker \rho=\{1\}$ therefore $F_r$ is the function field of the curve $X$.
Notice also that since $F_1$ is rational $F_0=F_r^{G_1(P)}\subseteq F_1$,
 is also rational.

Let $Q_i$ be the restriction at the intermediate fields $F_i$ of the place  $P$ of $F_r$.
 We will denote be $H_{F_i}(Q_i)$ the Weierstrass semigroup of the field $F_i$ at 
the place $Q_i$. 
We have assumed that  $F_1=k(f_1)$ is a rational function field and thus the 
semigroup of $P_1$ is $H_{F_1}(Q_1)=\Z_+$. This gives us that 
$|\ker \rho_1 | \Z_+ \subset H_{F_r}(P)$.

There are now two cases: Either the group $G_1(P)$ acts trivially on $f_1$ therefore 
$\ker\rho_1=G_1(P)$, or  
the group $G_1(P)$ might act on 
$f_1$ by translations (this is the only possible Galois action of a $p$-group on a rational 
function field see \cite{Val-Mad:80}). The function field generator of $F_0$
 {gives 
rise to  pole numbers that are multiplies of $|G_1(P)|$. Since $|\ker \rho_1|$ divides $ |G_1(P)|$
 all pole numbers coming from the function field generator of $F_0$ 
are already in $|\ker \rho_1 | \Z_+ \subset H_{F_r}(P) $.
Notice that $0$ is a jump in the representation filtration only if $f_1$ is not 
$G_1(P)$-invariant. In that case $F_0(f_1)=F_1$.

Suppose that $1$ is the first jump of the representation filtration, i.e. take $i_1=1$.
We proceed to the first $m_{i_1}$ such that $m_{i_1} \not\in |\ker\rho_1| \Z_+$,
 i.e. we choose the minimal $m_i$ with the property that $|\ker\rho_1| \nmid m_{i}$ 
for  $1< i$.
This means that the element $f_{i_1}$ is not $\ker \rho_1$-invariant  and 
$\ker \rho_1 \supsetneq \ker \rho_{i_1}$, i.e. the first { jump in the representation filtration} 
appears at the first $m_{i_1}$ that is not in $|\ker \rho_1| \Z_+$. The extension $F_{i_1}/F_1$
is elementary abelian  with Galois group $\ker\rho_1/\ker \rho_{i_1}$. 
The element $f_{i_1}$ 
 is not an element in $F_1$ since 
every element in $F_1$ is $\ker\rho_1$-invariant. 

Suppose now that $i_1\neq 1$. Then $\ker \rho_0=\ker \rho_1=\cdots=\ker \rho_{i_1}$, with $m_\nu \in |\ker \rho_{i_1}|\mathbb{Z}_+$ for all $\nu\leq i_1$. We choose $m_{i_2}$ minimally such that $m_{i_2}\not \in |\ker\rho_{i_1}| \Z_+$, i.e. $|\ker\rho_{i_1}| \nmid m_{i_2}$. This means that   the corresponding funtion  $f_{i_2}$ is not  $\ker \rho_{i_1}$-invariant and $\ker \rho_{i_1}\supsetneq \ker \rho_{i_2}$. Notice that $F_0=F_1=\cdots=F_{i_1}$.  The extension $F_{i_2}/F_{i_1}=F_{i_2}/F_1$
is elementary abelian  with Galois group $\ker\rho_{i_1}/\ker \rho_{i_2}$, but the element $f_{i_2}$ 
 is not an element in $F_{i_1}$ since  it is not $\ker\rho_{i_1}$-invariant. 

 The result is proved the same way going up step by step in terms of Artin-Schreier extensions.}
\end{proof}

\begin{remark}
 The extensions $F_{i_\nu+1}/F_{i_\nu}$ can be described as an Artin-Schreier
 extensions.
For example {} if $\nu=1$, with $i_1\neq0$, we have
 $\mathrm{Gal}(F_{i_1+1}/F_1)=\oplus_{h} \Z/p\Z$
for each cyclic direct summand of $\mathrm{Gal}(F_{i_1+1}/F_1)$. Notice that $F_0=F_1=\cdots=F_{i_1}$.
Then, for the action of $\sigma \in\oplus_{h} \Z/p\Z $ on $f_{i_1+1}$ we have:
\[
 \sigma^\ell (f_{i_1+1})=f_{i_1+1} + \sum_{0\leq \mu \leq i_1}   c_\mu (\sigma) f_{\mu},
\]
 and $c_\mu (\sigma)\in k$.
We have that $\sum_{0\leq \mu \leq i_1}    c_\mu(\sigma) f_{\mu}$ is $\ker \rho_{1}$-invariant. 
Thus, we set 
\begin{equation}\label{newprimes}
 Y_{i_1+1}=\frac{f_{i_1+1}}{\sum_{0\leq \mu \leq i_1}    c_\mu(\sigma) f_{\mu}}, 
\end{equation}
and observe that 
\[
 \sigma_\nu^\ell (Y_{i_1+1})=Y_{i_1+1}+\ell.
\]
This gives us that $Y_{i_1+1}$ is a set of Artin-Schreier generators for the extension $F_{i_1+1}/F_1$ and 
the elements $Y_{i_1+1}$ satisfy equations of the form $Y_{i_1+1} ^p - Y_{i_1+1}= a_{i_1}$ for some element 
in $a_{i_1}\in F_1=F_{i_1}$.  But actually we need only one extra generator $f_{i_1+1}$ in order to express
$F_{i_1+1}=k(f_1,f_{i_1+1})$.
\end{remark}

\begin{remark}

Observe that Artin-Schreier extension given  by eq. (\ref{newprimes})
introduces more ramification primes. Indeed, every zero of the quantity 
$\sum_{0\leq \mu \leq i_1}    c_\mu(\sigma) f_{\mu}$  is a ramified prime in the 
extension $F_{i_1+1}/F_{i_1}$. Also the valuation of $Y_{i_1+1}$ 
at the prime $P$  in question is {
$v_P(Y_{i_1+1})=v_P(f_{i_1+1})-v_P(\sum_{0\leq \mu \leq i_1}    c_\mu(\sigma) f_{\mu})$.}
\end{remark}

\begin{remark} \label{two-AS}
It seems natural to ask the following question: Let $L/K$ be a Galois 
extension with Galois group a   $p$-group.  Let $P$ be a  wild, totally ramified place 
of $L$. If a pole number in $H(P)$ is divisible by $p$, is it true 
that the corresponding element $f\in L$ is $\mathrm{Gal}(L/K)$-invariant?
The answer to this question is negative as one sees by the following example:
Consider the totally ramified  elementary abelian extension of the rational function field given 
by the two Artin-Schreier extensions: 
\[
 y_1^p-y_1=f_1(x),\qquad y_2^p-y_2=f_2(x),
\]
where $f_i(x) \in k[x]$ for $i=1,2$.
Consider the $\Z/p\Z$-extensions of $K_1:=k(x,y_1),K_2:=k(x,y_2)$ of the rational function field. 
In those extensions only the infinite place $\infty$ of $k(x)$ is ramified. Moreover it is known 
\cite[Korollar 1]{StiII} that 
the Weierstrass semigroups at the places of $K_i$ { that are above $\infty$}  are $p\Z_++m_i\Z_+$ respectively where 
$m_i=\deg f_i$.
We have the following tower of fields and corresponding semigroups:
\[
 \xymatrix{   k(x,y_1,y_2) \ar@{-}[d]_{\frac{\Z}{p\Z}} \ar@{-}[dr]^{\frac{\Z}{p\Z}}    & &  H(P)  \ar@{-}[d] \ar@{-}[dr] & \\
k(x,y_1) \ar@{-}[d]_{\frac{\Z}{p\Z}}  & k(x,y_2)  \ar@{-}[dl]^{\frac{\Z}{p\Z}}   &   p \Z_++ m_1 \Z_+  \ar@{-}[d]&
p\Z_+ + m_2\Z_+ \ar@{-}[dl] \\
 k(x) & &   \Z_+ & 
}
\]
Therefore,  the Weierstrass semigroup $H(P)$ of the place $P$ of $k(x,y_1,y_2)$ that is above $\infty$,
has $p^2\Z_+ + pm_1\Z_+ + pm_2 \Z_+$ as a subsemigroup. We have $\mathrm{div}_\infty(y_i)=pm_i$
but $y_2$ is not in $k(x,y_1)$, i.e., it is not $\mathrm{Gal}(k(x,y_1,y_2)/k(x,y_1))$-invariant. 
Suppose that $m_1<m_2$. Recently, \cite{QuiRe} and \cite{Sti:09} showed 
that the jumps of  the ramification filtration
occur at $m_1$ and at $m_1 +p(m_2-m_1)$.
\end{remark}

The situation changes for the case of cyclic extensions of the rational function field as will see 
in the next section.
\subsection{Cyclic totally ramified Galois {extensions of function fields}}
\label{cyc-section}

We will need the following:
\begin{proposition} \label{oooo}
 Let $L_m$ be a cyclic extension of the rational function field $L_0=k(x)$ of degree $p^m$. 
Assume also that every ramified place in the extension $L_m/L_0$ is ramified 
completely. 
Then 
\begin{enumerate}
 \item \label{11}
 there is a unique tower of intermediate fields
\[
 k(x)=L_0 \subset L_1 \subset \cdots \subset L_{m-1} \subset L_m
\]
such that $[L_j:L_{j-1}]=p$ and this extension is given in terms of an Artin-Schreier extension
\[
 L_j=L_{j-1}(y_j):\qquad y_j^p-y_j=B_j,
\]
where $B_j\in L_{j-1}$ has the following property (standard form):
\[
 \mathrm{div} B_j=A-\sum_{i} \lambda_{i,j} P_{i,j-1},
\]
with $A$ prime to the pole divisor of $B_j$ and $\lambda_{i,j}$ is either zero or a positive 
integer prime to the characteristic $p$.  Actually the integers $\lambda_{i,j}$ are equal to the 
valuations $-\lambda_{i,j}=v_{P_{i,j-1}}(B_j)=v_{P_{i,j}}(y_j)$.
In the above notation we are using the notation of \cite{Madden78}, namely in  
$P_{i,j}$ the index $j$ indicates that 
this prime is a prime of the field $L_j$, while the index $i$ shows  that lies above the $i$-th prime of
 $L_0$ that ramifies.
Moreover 
\item \label{12} The elements $y_1^{\alpha_1} y_2^{\alpha_2} \cdots y_n^{\alpha_m}$
 form a basis of $L_m$
over $L_0$, where $0\leq \alpha_i <p$. 
\item \label{13}  For every ramified place 
  $P$  in
 $L_m/L_0$  and integer $a\in \mathbb{N}$ we
 have the decomposition:
\[
 L(aP)=\bigoplus_{0\leq a_1,\ldots,a_m<p} \big( L(aP) \cap y_1^{\alpha_1} y_2^{\alpha_2} \cdots y_m^{\alpha_m} k(x) \big),
\]
where $y_1^{\alpha_1} y_2^{\alpha_2} \cdots y_m^{\alpha_m} k(x)$ denotes the 1-dimensional 
$k(x)$-vector space generated by $y_1^{\alpha_1} y_2^{\alpha_2} \cdots y_m^{\alpha_m}$.
\end{enumerate}
\end{proposition}
\begin{proof}
 Part (\ref{11}) is theorem 2 in \cite{Madden78}, part (\ref{12}) is lemma 3 in the same article and part (\ref{13})
follows by observing that the divisor $a^{-1}$ on page 311 can be replaced by any 
$\mathrm{Gal}(L_m/L_0)$-invariant 
divisor.  
\end{proof}
\begin{definition}
We will denote by 
$G(i)$ the subgroup of $G_1(P)$ of order $p^{m-i}$.
\end{definition}

\begin{proposition}\label{lambda_j}
Fix a place $P_{i_0,0}$ of the rational function field $L_0$
that ramifies {totally} in the extension $L_m/L_0$. 
 Then the set of places $P_{i_0,j}$ of 
$L_j$  that are above $P_{i_0,0}$ will be denoted by $P_{j}$ and 
 $\lambda_{j}=\lambda_{i_0,j}$,  $1\leq j\leq m$.

 Let $H_{L_j}(P_j)$ denote the Weierstrass semigroup at the place $P_j$ of the function field $L_j$. 
We have:
\begin{equation} \label{semigroupCyc}
 H_{L_j}(P_j) \subseteq  p^{j} \Z_+  + p^{j-1}  \lambda_{1} \Z_+ + \cdots +\lambda_{j} \Z_+, 
\end{equation}
In particular 
\[
 H_{L_m}(P_m) \subseteq  p^{m} \Z_+  + p^{m-1}  \lambda_{1} \Z_+ + \cdots +\lambda_{m} \Z_+
\]
If moreover $P_{i_0,0}$ is the only place that ramifies in the extension $L_m/L_0$ then 
the above inclusions are equalities.
\end{proposition}
\begin{proof}
Assume that  $\mu<j$ and let $y_\mu$ the generator of $L_\mu$ over $L_{\mu-1}$. 
Observe that on the function field $L_j$ we have
 $v_{P_j}(y_\mu)=-p^{j-\mu}\lambda_{\mu}$.
According to prop. \ref{oooo} part (\ref{13}) every element in $ \cup_{\nu \geq 1}L(\nu P_m)$ 
comes from an element of the form
 $y_1^{\alpha_1} y_2^{\alpha_2} \cdots y_m^{\alpha_m}f(x)$
and these valuations at $P_m$ form the semigroup 
 $p^{m} \Z_+  + p^{m-1}  \lambda_{1} \Z_+ + \cdots +\lambda_{m} \Z_+$.
Therefore $H_{L_m}(P_m)\subseteq  p^{m} \Z_+  + p^{m-1}  \lambda_{1} \Z_+ 
+ \cdots +\lambda_{m} \Z_+$.

We don't know if the elements $y_i$ have other poles so they may not 
be in $L(\kappa P_m)$ for some integer $\kappa$. But if there is only one ramified point 
then the pole divisors of $y_i$ is supported only at  $P_m$ and $\lambda_i p^{m-i}$ is in the 
semigroup.    
\end{proof}
Using the description of $L_m$ we are able to see the following
\begin{proposition} \label{29}
 If $\mathrm{Gal}(L_m/L_0)$ is cyclic of order $p^m$ and $H(P)$ is the Weierstrass semigroup 
at a totally ramified point { $P$ of $L_m$},  and $\mathrm{div}_\infty(f)=aP$ for an integer $a$ such that $p^\kappa \mid a$ 
then $f$ is in $L_{m-\kappa}$, i.e. $f$ is fixed by the unique subgroup of $\mathrm{Gal}(L_m/L_0)$, 
of order $p^\kappa$. 
\end{proposition}
\begin{proof}
 Using the valuation of elements of the form 
$y_1^{\alpha_1} y_2^{\alpha_2} \cdots y_m^{\alpha_m}$
we see that if  $p^\kappa \mid a$ then $f$ is a $L_0$ linear combination
of elements $  y_1^{\alpha_1} y_2^{\alpha_2} \cdots y_m^{\alpha_m}$ with {
$\alpha_m=\cdots=\alpha_{m-(\kappa-1)}=0$} that are in $L_{m-\kappa}$. 
\end{proof}

In eq. (\ref{semigroupCyc}) we introduced the sequence of numbers 
$\lambda_1 < \lambda_2< \ldots $ corresponding to the valuations 
of the generating elements $y_i$ of the field $L_m$. 
Is there any relation to the sequence $m_1< m_2< \ldots$ of pole numbers?
What is the relation with the jumps of the representation filtration?

Assume that the jumps of the representation filtration appear at the integers
$j_1,j_2,\ldots,j_n$, i.e. $\ker \rho_{j_i} > \ker \rho_{j_i+1} $ for all $i=1,\ldots, n.$ 
We have the following:
\[
 G_1(P)= \ker \rho_1=\ker \rho_2=\cdots=\ker \rho_{j_1} > \ker \rho_{j_1+1}=\cdots=
\ker\rho_r=\{1\}.
\]
Notice that the kernel of {$ \rho_r$ is trivial since $\rho_r$ is a faithful representation} 
by proposition \ref{fai-rep-lemma}. 
Observe that the  number of jumps in both the representation and the ramification filtrations 
 is the same, since in both filtrations the elementary 
abelian quotients at indices where a jump occur, are just cyclic groups of order $p$,
therefore  $m$, the number of jumps in the ramification filtration, equals $n$, the number of jumps in the representation filtration. 
The jumps at the representation filtration occur at 
$j_1,j_2,\ldots,j_{n}$, and $\ker \rho_{j_\nu}=G(\nu)$.

Let $j_\nu$ be an index such that $G(\nu)=\ker\rho_{j_\nu} \varsupsetneq \ker \rho_{j_\nu+1}=
G(\nu+1)$.
Since we are in the cyclic  group  case we have 
$\ker\rho_{j_\nu}/\ker \rho_{j_\nu+1}\cong \Z/p\Z$.
Let $\sigma$ be a generator of the later quotient. We can write the action on $f_{j_\nu+1}$ as 
{
\[
 \sigma f_{j_\nu+1}=f_{j_\nu+1} + \sum_{\mu\leq j_\nu} c_\mu(\sigma) f_\mu.
\]
}
{
Since all $f_\mu, \mu \leq j_\nu$ are $\sigma$-invariant
the element 
\begin{equation} \label{YYdef} 
 Y_{j_\nu+1}=\frac{f_{j_\nu+1}}{\sum_{\mu \leq j_\nu} c_\mu(\sigma) f_\mu}
\end{equation}
 is acted on by 
$\sigma$ as 
$
 \sigma Y_{j_\nu+1}=Y_{j_\nu+1}+1,
$
}
and it is an Artin-Schreier generator of the extension $F_{j_\nu+1}/F_{j_\nu}$. 
\begin{definition} \label{def:29}
Consider a tower $F =F_n \supseteq  F_{n-1} \supseteq  \cdots \supseteq  F_1 \supseteq 
  F_0$ of fields. 
On this tower of fields we consider the sequence of jumps
$F_n \varsupsetneq F_{j_m} \varsupsetneq F_{j_{m-1}} 
 \varsupsetneq \cdots \varsupsetneq F_{j_1} = \cdots= F_0$.
In the notation of proposition \ref{oooo} the fields $F_{j_\nu+1}$ are just the 
fields $L_\nu$. 
Fix a place $P$ of $F$ and let $Q_i$ be the restriction of $P$ on the field $F_i$.
Notice that the places $P_\nu$ of $L_\nu$ defined in proposition \ref{oooo} are equal to 
the places $Q_{j_{\nu}+1}$.  
 
An element $f_j \in F_j$ can be seen as an element in all fields $F_{j'}$ with 
$j' \geq j$. 
For $j\geq i$ we will denote by { $0>-m_{j,i}$  the valuation  of $f_i$ as an
 element in $F_j$ at the place $Q_j$, where the  function $f_i$ corresponds to the pole number 
$m_i$ of $H(P)$. We are interested in the elements $f_{j_\nu+1}$, i.e. for the indices 
where a jump occurs in the representation filtration. 
Notice that since $f_{j_\nu+1} \not \in \langle f_i: i\leq j_\nu \rangle$ the 
valuation $v_{Q_{j_\nu+1}}(f_{j_\nu+1})$ is prime to $p$ according to proposition \ref{29}.
Therefore the  number $m_{j_\nu+1,j_\nu+1}$ 
is prime to $p$  and we will denote it by $\mu_{\nu+1}$. }
\end{definition}

We compute using eq. (\ref{YYdef})
{
\begin{eqnarray} \label{valY}
v_{Q_{j_\nu+1}}(Y_{j_\nu+1}) & = &-m_{j_\nu+1,j_{\nu}+1}-
v_{Q_{j_\nu+1}}\left(\sum_{\mu\leq  j_\nu} c_\mu(\sigma) f_\mu\right) \nonumber \\
 &=& -\mu_{\nu+1}+\max \{ m_{j_\nu,\mu} : c_\mu(\sigma)\neq 0\}.
\end{eqnarray}
 Notice  that $ \lambda_{i,j}$ the valuations of generators $y_j$, 
of the fields given in the field tower of proposition \ref{oooo} part (\ref{11}),
 are  prime to the characteristic from the standart form hypothesis.
 Notice also that there are $\nu$ generators $Y_{j_\nu+1}$.
 We will abuse  the notation and denote with $Y_\nu$, $1\leq \nu \leq n$ the
  generators of the field extension $F_{j_\nu+1}/F_{j_\nu}$.} Now,     
since $Y_\nu$ and $y_\nu$ are both Artin-Schreier generators 
their valuation, { which is prime to the characteristic from the discussion above,} 
equals the conductor of the extension and must be equal, thus
the elements $Y_\nu$ and $y_\nu$ have the same valuation at 
$Q_{j_\nu+1}=P_\nu$. The absolute value of the valuation computed in eq. (\ref{valY}) equals 
$\lambda_\nu$, where $\lambda_\nu$ are the numbers introduced in proposition \ref{lambda_j}.
 
\begin{corollary}
We keep the notation from definition  \ref{def:29}.
 The valuations $\lambda_\nu$ of the elements $y_\nu$ at $P_{j_\nu}$ 
satisfy
$\lambda_\nu \leq \mu_\nu$. Equality holds only if there is only one place $P$ 
ramified in extension $F_n/F_0$, i.e., we are in the case of a Katz-Gabber cover with cyclic group of order $p^n$.
\end{corollary}
\begin{proof}
The inequallity is a direct consequence of eq. (\ref{valY}). 
If $\lambda_\nu=\mu_\nu$ then all coefficients $c_\mu(\sigma)$ $\mu>0$ in eq. (\ref{valY}) are zero, 
i.e. only the coefficient $c_0(\sigma)$ of the constant function $f_0$ can be non-zero.
Then $Y_\nu=f_{\nu}$ and there is only one place ramified in the extension $F_{\nu+1}/F_{\nu}$.
\end{proof}
\begin{proposition}
Let $P$ be as in definition \ref{def:29}.
The jumps of the ramification filtration at $P$ are exactly at the the integers $\lambda_\nu$,
 { where  $\lambda_\nu$ are the 
valuations of the elements $y_\nu$ at $P_{\nu}$.}
\end{proposition}
\begin{proof}
 Observe that  by \cite[proof. Prop. 2.3]{KontoZ}, the integer $\lambda_\nu$ equals the lower 
jump of the representation filtration $ \ker\rho_{j_\nu}/\ker \rho_{j_\nu+1}\cong \Z/p\Z$
 (keep in mind that in that paper  the enumeration is decreasing).
Notice that $ \ker\rho_{j_\nu}=G(\nu)$ and $\ker \rho_{j_\nu+1}=G(\nu+1)$.
The ramification filtration of $G(\nu)$ is given by 
$G(\nu)_\mu=G_\mu(P) \cap G(\nu)$.
Observe that the groups $G(\nu)$ are elements of the ramification filtration therefore there is no need 
to consider the upper ramification filtration, in order to relate the ramification filtrations of 
$G_1(P)$ and $G_1(P)/G(\nu)$ \cite[corollary after prop 3. IV.1]{SeL}, i.e.,
\[
 \left(G_1(P)/G(\nu)\right)_j=G_j(P)/G(\nu). 
\]
In the same way we have
\[
 \left(\frac{G(\nu)}{G(\nu+1)} \right)_j=\frac{G(\nu)_j}{G(\nu+1)}=\frac{G_j(P) \cap G(\nu)}{G(\nu+1)}.
\]
This implies that 
$\lambda_\nu$ are the jumps of the ramification filtration.
\end{proof}

\begin{example}
We consider now the case of a a cyclic group extension $L/L_0$
 of order $p^n$ with only one full ramified 
place such that only $L_0$ is rational. 

 We have the following tower of $\Z/p\Z$-cyclic extensions:
\[
\xymatrix{
 L \ar@{.}[d]   & P \ar@{.}[d]   &   p^n \Z_+ +  p^{n-1} \lambda_1 \Z_+ + p^{n-2} \lambda_2 \Z_+ +\cdots + \lambda_n \Z_+
 \ar@{^{(}.>}[d]   \\
L_2=L^{G(2)} \ar@{-}[d]_{\frac{\Z}{p\Z}} & P_2  \ar@{-}[d]&  p^2 \Z_+ +p \lambda_1\Z_+ + \lambda_2 \Z_+ \ar@{^{(}->}[d]\\
L_1=L^{G(1)} \ar@{-}[d]_{\frac{\Z}{p\Z}} & P_1  \ar@{-}[d]&  p\Z_+  +\lambda_1 \Z_+ \ar@{^{(}->}[d]\\
L_0=L^{G_1(P)} & P_0  &\Z_+
}
\]
The place  $P$ is fully ramified in extension $L/L^{G_1(P)}$ and we will denote by $P_i$ 
the corresponding place of  $L^{G(i)}$.
The field $L^{G_1(P)}$ is assumed to be rational and has Weierstrass semigroup $\Z_+$. 
The field $L^{G(1)}$ is an Artin-Schreier extension of $L^{G_1(P)}$ 
and  the Weiestrass semigroup  at $P_1$ has a part $p \Z_+$ coming as the 
semigroup of $\Z/p\Z=G_1(P)/G(1)$-invariant elements plus a new non invariant element 
that is non divisible by $p$, the $\lambda_1 \Z_+$ semigroup.
 The Weierstrass semigroup at $P_1$  is 
$p \Z_+ + \lambda_1 \Z_+$. 

The next step is to consider the Weierstrass semigroup of $P_2$. It has a part 
$p^2 \Z_+ + p \lambda_1 \Z_+$ coming from the $G(1)/G(2) \cong \Z/p\Z$-invariant
 elements and also 
some extra elements that should contain a semigroup of the form 
$ \lambda_2 \Z_+$, where $(\lambda_2,p)=1$. 
But  eq. (\ref{semigroupCyc}) gives us that 
$H(P_2) \subseteq p^2 \Z_+ + p \lambda_1 \Z_+ + \lambda_2\Z_+$ and since
there is a unique ramified place the above inclusion is actually an equality.

This way we can go all the way up to $L$ 
and find that the Weierstrass semigroup  is given by 
\[
 H(P)=p^n \Z_+ +  p^{n-1} \lambda_1 \Z_+ + p^{n-2} \lambda_2 \Z_+ +\cdots + \lambda_n \Z_+.
 \]
%
\end{example}
\subsection{Holomorphic differentials in the cyclic  totally ramified
 case}
\label{hol-diff}

In this section we are going to use known bases for the 
space of holomorphic differentials in order to express the gaps of the Weierstrass 
semigroups in terms of the Boseck invariants of the curves. In  \cite{karan}  we have defined Boseck invariants to be the values  $\Gamma_k(m)$ in definition \ref{cyclic} below 
and we expressed the Galois module structure of holomorphic $m$-polydifferentials in terms of them.
Here we we omit the $m$ from the notation since we are interested only for the $m=1$ case. 
The motivation for their name was the work of Boseck \cite{boseck}
and in what follows we will call them the {\em Boseck invariants} of the curve.

We assume that $F$ is a cyclic, totally ramified extension of order $p^n$ of the rational function field.

Let $k$ be an integer with $p$-adic expansion 
\begin{equation}\label{k}
 k=a_1^{(k)}+a_2^{(k)}p+\cdots+ a_{n}^{(k)}p^{n-1}.
\end{equation}
The set $w_k=y_1^{a_1^{(k)}}y_2^{a_2^{(k)}}\cdots y_n^{a_n^{(k)}}$, 
$0\leq k \leq p^n-1$ is an  $k(x)$-basis $E$-basis of $F$ 
 \cite{vm}.
The valuations of the basis elements $w_k$ are given by 
\begin{equation}\label{basis}
v_{P_{i,n}}(w_k)=-\sum_{j=1}^n a_j^{(k)} \lambda_{i,j}  p^{n-j},
\end{equation}
{ where the $\lambda_{i,j}$ are given in proposition \ref{oooo}, part (\ref{11}).}


We denote the ramified places of $k(x)$, with $(x-\alpha_i),\ 1\leq i\leq s$, since in a rational function field  every ramified place 
corresponds to an irreducible polynomial, which is linear since the 
the field $k$ is algebraically closed.   We set
\begin{eqnarray*}
 g_k (x)=\prod_{i=1}^{s} (x-\alpha_i)^{\nu_{ik}(1)}.
\end{eqnarray*}
\begin{definition}\label{cyclic}
For $k=0,1,\ldots,p^n -1,$ we define
\[
\Gamma_k:=\sum_{i=1}^s \nu_{ik},
\]
where
\begin{equation}\label{niks}
\nu_{ik}=\biggl\lfloor \frac{\delta_i - \sum_{j=1}^{n}a_j ^{(k)} \lambda_{i,j}p^{n-j}}{p^{n}}\biggr\rfloor,
\end{equation}
where   $\delta_i$ denotes  the different exponent $\delta(P_{i,n,\mu}/(x-\alpha_i))$ and the $1\leq j \leq n$ runs over the intermediate fields in 
 the tower of proposition \ref{oooo}, part \ref{11}. Finally 
\begin{equation}\label{rho}
\rho_{i}^{(k)}=\bigl(\delta_{i}-\sum_{j=1}^n a_j^{(k)} \lambda_{i,j} p^{n-j}\bigr)-\nu_{ik} \cdot p^{n},
\end{equation}
is the remainder of the division of the quantity $\delta_i +v_{P_{i,n}}(w_k)$ by $p^n$.
\end{definition}

\begin{proposition} \label{holbasis-cyc} 
Let $X$ be a cyclic extension of degree $p^n$ of the rational function field.
The set 
\[
\left\{
\omega_{k\nu}^{(\alpha_i)} =(x-\alpha_i)^{\nu^{(k)}} g_{k}(x)^{-1}w_{k}dx: 0 \leq \nu^{(k)} \leq \Gamma_k -2, 0\leq k \leq p^n-2\right\}
\]
forms a basis for the set of holomorphic differentials for a cyclic extension of the 
rational function field of order $p^n$.
\end{proposition}
\begin{proof}
 We take the basis of \cite[Lemma 10]{karan}, set $m=1$ and modify it in order to evaluate
holomorphic differentials in the ramified primes of the extension.
The same construction is given by Garcia in \cite[Theorem 2, Claim]{garciaelab} 
where  the  elementary abelian, totally ramified case is studied. 
The proof is identical to the one given there.
\end{proof}

\begin{remark}\label{maxel}
Observe from eq. (\ref{niks}) that if $\Gamma_k=0$, then $\nu_{ik}=0$ 
for all $i=1,\ldots,r$  and that could happen only when $k =p^n -1$ since $\nu_{ik} \geq 0$. Thus 
 $\Gamma_k$'s attain the maximum value when $k=0$, i.e $\Gamma_0 \geq \Gamma_k$, for 
 all $0\leq k \leq p^n-1$ and  for  $0\leq k \leq p^n-2$,  we have that $\rho_i^{(k)}\leq p^n-2$.
\end{remark} 

We will need the following definition. 
For more information concerning the theory of linear series 
the reader is referred to \cite{Torres00}.
\begin{definition}
 Let $\mathcal{D}$ be a linear series of degree $d$ and dimension $r$, i.e.,
a linear subspace of $\mathbb{P}(L(E))$. The linear series can be seen as 
a set of linear equivalent divisors.
We can form a decreasing sequence 
\[
 \mathcal{D}_i(P)=\{D \in \mathcal{D}:  D \geq i P\}.
\]
The {\em sequence of $(\mathcal{D},P)$-orders}  is  a sequence of integers:
\[
 j_0^{\mathcal{D}}(P)<\cdots <j_r^{\mathcal{D}}(P),
\]
 such that $\mathcal{D}_j(P) \varsupsetneq \mathcal{D}_{j+1}(P)$ for all $j$ in the above sequence. 

For all but finite points of $X$ the $(\mathcal{D},P)$-orders do not depend on $P$. The 
exceptional points are called $(\mathcal{D},P)$-Weierstrass points. 

Let $\mathcal{E}_{\mathcal{D}}$ be the generic $(\mathcal{D},P)$-order sequence. 
We will call it the {\em order sequence of the linear series}.

A curve will be called {\em classical} with respect to the linear series $\mathcal{D}$ 
if and only if $\mathcal{E}_{\mathcal{D}}=\{0,\ldots,r \}$.
\end{definition}

We have the following:
\begin{proposition}\label{polenumbers}
Let $X$ be a cyclic extension of degree $p^n$ of the rational function field.
The gap sequence at the ramified primes $P_i$ that lies over the place $(x-\alpha_i)$, with $0\leq i \leq s$
of $k(x)$ is given by 
\[
\mathcal{G}(P_i)=\left\{\nu^{(k)} \cdot  p^n +\rho_i^{(k)}+1 |\ 0\leq k \leq p^n-2,\ 0\leq \nu^{(k)} \leq \Gamma_k  -2 \right\}, 
\]
where $\Gamma_k$ is the Boseck invariant associated to the extension $F/k(x)$.
\end{proposition}
\begin{proof}
This computation is based on the valuation of  the elements in the basis of holomorphic 
differentials given in proposition \ref{holbasis-cyc}.
Let $\rho_i^{(k)}$ be as in Eq. (\ref{rho}).
We denote  by   $j^K (P_i)=j^K$ the $(K,P_i)$-orders and $K$ the canonical linear series.
By computing the valuations at $P_i$ of the basis elements given in proposition \ref{holbasis-cyc},
we are able to compute that
the  $(K,P_i)$ orders are given by 
\begin{equation}
 \label{differentValues}
\left\{\nu^{(k)}\cdot  p^n +\rho_i^{(k)} |\ 0\leq k \leq p^n-2,\ 0\leq \nu^{(k)} \leq \Gamma_k  -2\right\}.
\end{equation}
In the above computation it is essential to  notice that for different values of 
$k, \nu^{(k)}$ $0\leq k \leq p^n-2, \ 0\leq \nu^{(k)} \leq \Gamma_k  -2$
the values $\nu^{(k)} \cdot  p^n +\rho_i^{(k)}$ are different, therefore the 
valuation of a linear combination of the differentials in proposition \ref{holbasis-cyc},  is just the minimal of the valuation of each summand.
We know that knowledge of  the $(K,P_i)$'s orders,  is equivalent 
to the knowledge of the $(0,P_i)$'s gaps, i.e. the gap sequence at
$P_i$ will be given by $j^K_{\mu}(P)+1$. 

 Notice that every element has 
a different valuation with respect to the place $P_i$. Notice also that the dimension 
of the holomorphic differentials is $g$ the same number as the number of the  gaps. 
Also the $n=1$ case was studied by Boseck in \cite[Satz. 17]{boseck}
\end{proof}

\begin{remark}
Boseck in his seminal paper \cite[Satz 18]{boseck},  where the $n=1$ case is studied, states that as  
as $k$ takes all the values $0\leq k \leq p-2$ the remainder of the Boseck's basis 
construction $\rho_i^{(k)}$ takes all the 
values $0\leq \rho_i^{(k)} \leq p -2$  and thus all the numbers $1,\ldots,p-1$ are gaps.
 This is not entirely correct as we will show in example 
\ref{31acc}. The problem appears if there  is exactly  one ramified place in the Galois 
extension.  
\end{remark}

\begin{lemma} \label{smallGaps}
 If all $\Gamma_k \geq 2$ then all numbers  $1,\ldots,p^n-1$ are gaps.
If  there are  Boseck invariants $\Gamma_k=1$, then the set of gaps smaller than 
$p^n$ is exaclty the set $\{\rho_i^{(k)}+1|\ k:0 \leq k \leq p^n-2, \Gamma_k \geq 2\}$.
\end{lemma}
\begin{proof}
Recall that 
in eq. (\ref{rho}) the elements $\rho_i^{(k)}$ were defined to be  the remainder of the division of 
$\delta_i -\sum_{j=1}^n a_j^{(k)} \lambda_{i,j} p^{n-j} $ by $p^n$. 
As $k$ runs in $0\leq k \leq p^n-2$ the $\rho_i^{(k)}$ run in $0,\ldots,p^n-2$.
Indeed, let us define the function
\[
\Psi:
 \{0,\ldots,p^n-2\} \rightarrow \{0,\ldots,p^n-2\},
\]
\[
 k \mapsto \rho_i^{(k)}.
\]
Since the expressions in eq. (\ref{differentValues}) are all different the function 
$\Psi$ is onto.

But the $\Gamma_k$ that are equal to $1$ have to be excluded since they give not rise 
to a holomorphic differentials in proposition \ref{holbasis-cyc}, see \cite[Eq. (21)]{karan} 
and  example \ref{31acc}.
\end{proof}
\begin{remark}
Notice that elements $\Gamma_k=1$ can appear only for primes $p\geq \lambda_{i,j}$ and 
only if there is only one ramified place. But in our representation theoretic viewpoint 
the interesting case is the one of small primes. 
\end{remark}

\begin{example} \label{31acc}  We consider the now the case of an Artin-Schreier extension of the 
function field $k(x)$, of the form $y^p-y=1/x^m$. 
In this extension only the place $(x-0)$ is ramified with different exponent $\delta_1=(m+1)(p-1)$. 
The Boseck invariants in this case are 
\[
 \Gamma_k=\lf \frac{(m+1)(p-1)-km}{p} \rf \;\; \mbox{ for } k=0,\ldots,p-2.
\]
The Weierstrass semigroup is known \cite{StiII} to be  
$
m \Z_+ + p \Z_+$. 
Let us now try to find the small gaps by using lemma \ref{smallGaps}.
If $p<m$ then all numbers $1,\ldots,p-1$ are gaps. If $p > m$ then $m$ is a pole number smaller 
than $p$.  
Indeed, $\Gamma_{p-2}=1$ and the remainder of the division of 
$(m+1)(p-1)-(p-2)m$ by $p$ is $\rho^{(p-2)}=m-1$. But then 
$\rho^{(p-2)}+1=m$ is not a gap.  
\end{example}
\begin{definition} \label{symetric-semi}
In the literature, see Oliveira \cite{Oliveira1991}, semigroups $H$ of genus $g$ such the $g$-th gap is 
 $2g-1$ are called {\em symmetric}. 
\end{definition}\label{symetric-semi-pro}

A semigroup $H$ is called symmetric, because the symmetry is expressed in the semigroup in the following way
\begin{equation} \label{15aaaa}
n\in H \mbox{ if and only if }  2g-1-n \not \in H.
\end{equation}

\begin{lemma} \label{cyc-symmetric}
Let $X$ be a cyclic $p$-group extension of the rational function field with only one totally ramified point  $P_{i_0}$.
The bigest gap is equal to $(\Gamma_0-2)p^n+\rho^{(0)}+1=\delta -2 p^n+1$ and 
the { Weierstrass} semigroup at that point is symmetric. 
\end{lemma}
\begin{proof}
 Since there is only one ramified place $\nu_{i_0k}=\Gamma_k$.
 From remark \ref{maxel} and proposition \ref{polenumbers} we see that
the biggest gap is equal to $(\Gamma_0-2)p^n+\rho^{(0)}+1=\delta -2 p^n+1$, 
where $\delta$ is the  different exponent at the unique ramified place and is equal to the 
 degree of the different. 
For the last equality keep in mind that $\Gamma_0=\lf \frac{\delta}{p^n} \rf$ and 
$\delta=p^n \Gamma_0 +\rho^{(0)}$. 

A direct computation with Riemann--Hurwitz formula, shows that
\[
 2g-1 =-2p^n+1 +\delta,
\]
hence the semigroup is symmetric. 
\end{proof}
 Does   lemma  \ref{cyc-symmetric} holds for the general case when we have
  more than one totally
 ramified primes?
The answer is no, in general, but the following lemma  gives a  necessary and sufficient
 condition for this.
\begin{lemma}\label{cyc-symms}
 Assume that $X$ is as in lemma \ref{cyc-symmetric} but now there are $s\geq 2$ totaly 
ramified places $P_{i}$, with $1\leq i \leq s$. Let $\delta_i$ denote the 
different exponent at the prime $i$. 
 The semigroup at {some} ramified place $P_{i}$  is symmetric if and 
only if 
\begin{equation}\label{11wwee}
 \sum_{i'\neq i}\delta_{i'}=p^n\sum_{i'\neq i} \lf \frac{\delta_{i'}}{p^n} \rf. 
\end{equation}
Moreover, let $\lambda_{i,j}$ 
be the  valuations of the genarators of 
the field extensions given in Proposition \ref{oooo} part (\ref{11}). 
Then eq. (\ref{11wwee}) holds if and only if
  $\lambda_{i'} \equiv -1 \mod p^n$ for all $i'\neq i$. 
\end{lemma}
\begin{proof}
  Fix an $i$. Then from remark \ref{maxel} 
and proposition \ref{polenumbers}, we see that
the biggest gap at $P_i$ is $(\Gamma_0-2)p^n+\rho^{(0)}_i+1$ and using one more time the Riemann--Hurwitz formula,
 we see that this gap equals to $2g-1$ is equivalent to the condition
 \begin{equation}\label{gcase}
  (\Gamma_0-2)p^n+\rho_{i}^{(0)}+1=\sum_{i'=1}^s  \delta_{i'} -2 p^n+1.
 \end{equation}
Now since $\Gamma_0 =\sum_{i'=1}^s \lf \frac{\delta_{i'}}{p^n} \rf$, 
the left hand of eq. (\ref{gcase})  equals to
\[
p^n \sum_{i'\neq i} \lf \frac{\delta_{i'}}{p^n} \rf +\delta_{i}-2p^n +1.
\]
Thus in order eq. (\ref{gcase}) to be valid we should have
\begin{equation}\label{deltamod}
 \sum_{i'\neq i}\delta_{i'}=p^n\sum_{i'\neq i} \lf \frac{\delta_{i'}}{p^n} \rf. 
\end{equation}
The right hand of eq. \ref{deltamod} equals to $\sum_{i\neq i'}\bigl(\delta_{i'} -\rho^{(0)}_{i'}\bigr)$. Since $\rho^{(0)}_{i'}$'s are by definition non negatives, eq. \ref{deltamod} holds if and only if 
\[
\rho^{(0)}_{i'}=0 \Leftrightarrow \delta_{i'} \equiv 0 \mod p^n, \textrm{for every } i'\neq i.
\]
For $n=1$, $\delta_{i'} =(\lambda_{i',j} +1)(p-1)$ and the above condition is equivalent  to $\lambda_{i',j}\equiv -1 \mod p$, for $i'\neq i$.

For $n>1$, $\delta_{i'}=(p-1)\sum_{j=1}^n (\lambda_{i',j} +1) p^{n-j}$ (see \cite[p.110]{vm}), and the above condition is equivalent to
$\sum_{j=1}^{n-1} (\lambda_{i',j} +1) p^{n-j}\equiv \sum_{j=1}^n (\lambda_{i',j} +1) p^{n-j} \mod p^n$, for $i'\neq i$, or
$\lambda_{i',j}\equiv -1 \mod p^n$, for $i'\neq i$.
%
  
%
%


%
\end{proof}

\begin{corollary}
 Let $X$ be as in lemmata  \ref{cyc-symmetric},\ref{cyc-symms},
 i.e. assume that there exists a $P_{i_0} $ that is totally ramified and 
that its Weierstrass semigroup is symmetric. Then the Weierstrass sequence up to $2g$, at
  the  place $P_{i_0}$ that lies over the place $(x-\alpha_{i_0})$, 
 is given by
 \[
  H(P_{i_0})=\{2g-1 -a|a \in \mathcal{G}(P_{i_0})\}, 
 \]
where $\mathcal{G}(P_{i_0})$ is the gap sequence at $P_{i_0}$ from proposition \ref{polenumbers}.
\end{corollary}
\begin{proof}
 This is a direct consequence of lemmata  \ref{cyc-symmetric}, \ref{cyc-symms} and 
eq. (\ref{15aaaa}).
\end{proof}

\begin{remark} \label{Morrison}
In the theory of numerical semigroups the following construction is frequently used in order to describe 
the semigroup \cite{MorrisonPinkham}: Let $d(P)$ be the least positive element of $H(P)$. All 
elements $\mu \in \{1,\ldots,d(P)-1\}$ are gaps and for every such $\mu$ we denote by $b_\mu(P)$ 
the minimal 
element of $H(P)$ such that $b_\mu (P)\equiv \mu \mod \ d(P)$. 
This means that $b_\mu(P)=\nu_\mu(P)d(P)+\mu$, and $\nu_\mu(P)=
\lf \frac{b_\mu(P)}{d(P)} \rf$ equals the number of gaps that 
are congruent to $\mu$ modulo $d(P)$.  

Assume that the smallest pole number is  $p^n$ and that there is only one ramified place 
in the field $F_n/F_0$. 
Then the  integers $\nu_\mu(P)$ in this description of the semigroup are equal to the Boseck invariants
$\nu_\mu(P)=\Gamma_{\Psi^{-1}(\mu)}-1$ since both integers count the number of gaps that 
are equal to $\rho_i^{(k)}+1 \mod p^n$ by proposition \ref{polenumbers}.
%
Notice also that a theorem due to Lewittes \cite[th. 1.3]{MorrisonPinkham}, \cite[th. 5]{Lewittes63}
in characteristic zero, has an interpretation for the trivial group (the prime to $p$-part of a 
$p$-group acting on our curve):
\[
 g=\sum_{\mu=1}^{p^n-1} \nu_\mu=\sum_{\mu=1}^{p^n-1} (\Gamma_{\mu-1}-1),
\]
since $g$ is the trace of the trivial representation on holomorphic differentials 
and $\Gamma_\mu-1$ counts the 
number of gaps that are equivalent to $\mu$ modulo $p^n$. 
 This is equivalent to the formula proved in \cite[rem. 7]{karan} for $m=1$. 
\end{remark}

\begin{proposition}
We use the notation of theorem \ref{Lewittes}. Consider an Artin-Schreier cyclic extension
of the rational function field $F_0$, i.e.,  the Galois group is isomorphic to $\Z/p\Z$.  
Then for every $i$ there is an integer $k(i)$  with $0\leq i, k(i)\leq p-1$, such that 
$\lf\frac{mi}{p} \rf=\Gamma_{k(i)}-1$ for some Boseck invariant. 
\end{proposition}
\begin{proof}
According to theorem \ref{Lewittes} 
 the number of gaps equivalent to $mi$ modulo $p$ equals $\lf \frac{mi}{p}\rf$. 
By remark \ref{Morrison}
  this number equals to $\Gamma_k-1$ for the $k$ such that 
$\rho^{(k)}+1=mi \mod p$.
\end{proof}
\section{Semigroups}
\label{sec-semigroups}

Let $H\neq \mathbb{Z}_+$ be a semigroup of natural numbers and suppose that there is a natural number 
$n$  such that for all $s\geq n$ we have  $s\in H$. We select the minimal such number $n$, i.e. 
$n-1 \not\in H$. Observe that  when that semigroup is the Weierstrass semigroup of a curve,
 then the genus of the curve and the genus of the semigroup 
(i.e. the number of gaps of the semigroup) coincide, and we should have that 
 $g+1 \leq n \leq 2g$.  Let $d_1$ be an element  in $H$ and 
let $d_2$ be the minimal element in $H \backslash d_1\mathbb{Z}_+$ such that 
$(d_1,d_2)=1$.  Write $d_2=d_1 \lf \frac{d_2}{d_1} \rf  +u$, $0< u < d_1$.

\begin{lemma} \label{llee}
The numbers $d_1$, $d_2$ generate a subsemigroup 
$d_1 \mathbb{Z}_+  +d_2 \mathbb{Z}_+ \subset H$, such that every number 
$s>(d_1-1)(d_2-1)$ is 
in $d_1 \mathbb{Z}_+  +d_2 \mathbb{Z}_+ \subset H$.  In particular 
\begin{equation}
 \label{bound-n}
n\leq (d_1-1)(d_2-1). 
\end{equation}
If the numbers $d_1,d_2$ generate the semigroup $H$ then 
\begin{equation}\label{equal-n}
 n=(d_1-1)(d_2-1).
\end{equation}
\end{lemma}
\begin{proof}
Write $d_2=d_1\lf \frac{d_2}{d_1} \rf+ u$, $(u,d_1)=1$, $0<u<d_1$. In what follows we will describe the
 integers which  can be written as linear combination with $\Z_+$-coefficients of $d_1,d_2$. 

By considering elements of the form $\nu d_1+d_2$ we obtain all elements 
of the form $\left(\nu+\lf \frac{d_2}{d_1} \rf\right) d_1 +u$, $\nu\geq 0$.
Let $I_\lambda$ denote the interval $[\lambda d_1 , (\lambda+1)d_1)$. 
For $\lambda \geq \lf \frac{d_2}{d_1} \rf$ every interval $I_\lambda$
contains an integer that is equivalent to $u$ modulo $d_1$.

 By considering elements of the form 
$\nu d_1 + 2d_2$ we obtain all elements of the form 
$\left(\nu+2\lf \frac{d_2}{d_1} \rf \right) d_1 +  2 u$ etc.

In order to obtain an element with arbitrary residue modulo $d_1$ we have to consider 
all elements of the form $\nu d_1 + \mu d_2$, where $\mu$ takes all values in 
$0,\ldots,d-1$. Fix such an $\mu$. We consider  the combination
 \begin{equation} \label{mmm}
\left(\nu+\mu \lf \frac{d_2}{d_1} \rf\right) d_1 +\mu u
=  \left(\nu+\mu \lf \frac{d_2}{d_1} \rf  + \lf \frac{\mu u}{d_1} \rf \right) d_1 +
\left(\mu u- \lf \frac{\mu u}{d_1}\rf d_1\right).
 \end{equation}
This proves that for every $\lambda \geq 
 \left(\mu \lf \frac{d_2}{d_1} \rf  + \lf \frac{\mu u}{d_1} \rf \right)$
every interval $I_\lambda$ contains an integer that is equivalent to $\mu u$ modulo $d_1$.

The greatest value for $\mu$ is $d_1-1$ so for 
\[
 \lambda \geq (d_1-1) \lf \frac{d_2}{d_1} \rf  + \lf \frac{(d_1-1) u}{d_1} \rf =
(d_1-1) \lf \frac{d_2}{d_1} \rf  +u - 1. 
\]
For this value of $\mu$ 
the coefficient in front of $d_1$ in eq. (\ref{mmm}) is 
$(d_1-1)\lf \frac{d_2}{d_1}  \rf + u-1$, and this means that every natural number
 $\geq  d_1(d_1-1)\lf \frac{d_2}{d_1} \rf+ d_1(u-1)$ is in $H$.
But now we 
replace $d_1  \lf \frac{d_2}{d_1} \rf$ by $d_2 -u$
and we verify that:
\begin{equation} \label{11eee}
 d_1(d_1-1)\lf \frac{d_2}{d_1} \rf+ d_1(u-1)=(d_1-1)(d_2-1)+ u-1.
\end{equation}
We now observe that 
in the interval $I_{d_1-2}$ there is exactly one gap that is equivalent to  $d_1-u \mod d_1$.
The value of this gap equals 
$
 (d_1-1)(d_2-1)+ u-1-u=(d_1-1)(d_2-1)-1.
$
Therefore, all numbers greater than this gap are in $d_1 \Z_+ + d_2 \Z_+$. This means that 
\[
s  \geq  (d_1-1)(d_2-1) \Rightarrow s \in d_1 \Z_++d_2 \Z_2.
\]
Notice that we have proved that $ (d_1-1)(d_2-1)-1 \not\in d_1 \Z_+ + d_2 \Z_+$. 
If $H=d_1 \Z_+ + d_2 \Z_+$ then $(d_1-1)(d_2-1)-1$ is a gap and 
$n=(d_1-1)(d_2-1)$.
\end{proof}

Let $m_i$ be the sequence that enumerates the Weierstrass semigroup, $m_0=0$ is always 
a pole number. 
Observe that if $2g-1 \not \in H(P)$ then $n=2g$. 
Indeed, we  know that if the function field $F$ is not hyperelliptic then 
$m_i\geq 2i+1$ for $i=1,\ldots,g-2$ and $m_{g-1}\geq 2g-2$ \cite[lemma 1.25]{Torres00}.
This means that we have two cases for $m_{g-1}$, namely either 
$m_{g-1}=2g-2$ or $m_{g-1}=2g-1$. 

Let $K$ be the canonical linear series. 
Observe that at a generic point $P$ of a $K$--classical curve $X$, we have
$m_i (P) =g+i$, for  $i\geq 1$. In that case the gaps $\mathcal{G}(P)$ and the generic order sequence 
$\mathcal{E}(P)$ are classical and  they are equal to
\[
\mathcal{G}(P)=\{1,\ldots ,g\} \mbox{ and } \mathcal{E} (P)=\{0,\ldots,g-1\}.
\]
In what follows assume that $m_{g-1} =2g-2$ at  a point $P $ of the curve.  Notice 
 that this condtition  implies that the maximum gap at $P$  equals to $2g-1$.  Therefore this  leads to the study 
where the curve must satisfy \textit{at least one} of the following conditions: 
\begin{enumerate}   
\item the curve is not  $K$-classical,
\item   $P$ is a Weierstrass point.
\end{enumerate} 

In fact there is more to say. A more carefull  analisys indicates that if  $m_{g-1} =2g-2$, then the second condition, i.e.
 $P$ should  be a Weierstrass point, is always satisfied. Indeed, we distinguish the following cases: 

\textbf{Case1.}: The curve $X$ is not $K$-classical and $P$ is ordinary, meaning that
$2g-1 \in \mathcal{G}(P)\neq \{1,\ldots ,g\}$ and
\[
 \mathcal{G}(P)=\{\epsilon_0^K +1 , \ldots, \epsilon^K _{g-1}+1\}=\{j_{i}(P)^K+1|\ 0\leq i \leq g-1\}.
\]
Thus, we must have that $\epsilon_{g-1}^K = 2g-2 =j_{g-1}^K(P) $; 
This case cannot occur, see \cite[Lemma 2.31 p.30.]{Torres00}, \cite{Garcia:93}.

\textbf{Case 2.}: The point $P$ is a Weierstrass point with respect to K, and $X$ is a $K$-classical curve.
 That is $2g-1 \in \mathcal{G}(P) \neq \{1,\ldots ,g\}$, 
\[
 \mathcal{G}(P)=\{j_{i}^K(P)+1| \ 0\leq i \leq g-1\}\neq \{\epsilon_{i}^K+1|\ 0\leq i \leq g-1\},
\]
where $j_{g-1}^K(P) = 2g-2$ but $\epsilon_{g-1}^K=g-1 \neq  j_{g-1}^K(P) $.

\textbf{Case 3.}: The point $P$ is a Weierstrass point with respect to K and the curve $X$ is not K--classical. 
That is $g-1<\epsilon_{g-1}^K \lneq  j_{g-1}^K(P) $.

We have proved the following:
\begin{proposition}\label{genera}

If a curve has symmetric Weierstrass semigroup at a point $P$, then this point is a Weierstrass point.
Moreover if the Weierstrass semigroup at this point is generated by two elements, $d_1, d_2$, i.e.
 when we have  equality in lemma \ref{llee}, eq. (\ref{equal-n}), the genus of the curve is given by 
\[
g=\frac{(d_1 -1)(d_2-1)}{2}.
\]
If  the Weierstrass semigroup at this point cannot be generated by two elements, then 
\[
g<\frac{(d_1 -1)(d_2-1)}{2}.
\]

\end{proposition}

%

\begin{example}
%

Consider the case of  Artin--Schreier curve defined in example \ref{31acc} or the 
most general case
\[
y^q-y=f,
 \]
(where $f$ is a polynomial of degree  $m$ witch is prime to $p$, and $q=p^h$)
have genus $g=\frac{(q-1)(m-1)}{2}$.
Therefore,  using the theory of Boseck invariants we see that the biggest gap is 
\[
 \delta-2q+1=(m+1)(q-1)-2q+1=(m-1)(q-1)-1=2g-1,
\]
and  $n=2g$.

To the same conclusion we arrive using lemma \ref{llee}. 
Indeed, it is known \cite{StiII} that the Weierstrass semigroup is 
generated by $m,q$. Thus $n=(q-1)(m-1)=2g$. 
 
{Observe that certain function fields like the Hermitian 
function fields and their quotients $y^{q}+y=x^m$, $m \mid q+1$, $q=p^h$, satisfy $2n=g$. Also
 Matignon-Lehr curves  given by \cite[4.1]{CL-MM}, 
are given as Artin-Schreier extensions of the rational function field 
and satisfy $2n=g$. Moreover if $k=\mathbb{F}_{q^2}$, where $\mathbb{F}_{q^2}$ is a finite field with $q^2$ elements, then
Hermitian and their  quotients $y^{q}+y=x^m$, $m \mid q+1$, are certain maximal curves and can be viewed as the Artin-Schreier curves 
$y^q -y =f(x)$, where $f(x) \in \mathbb{F}_{q^2}[x]$ and $\mathrm{g.c.d.}(\deg f, p)=1$, see \cite[Theorem 5.4]{GarciaTafazolian}.
}


\end{example}

Here it is nice to point out the connection of maximal curves with Weierstrass semigroups:
Assume that $X$ is a maximal curve over $\mathbb{F}_{q^2}$ of genus $g$. Let 
$X(\mathbb{F}_{q^2})$ be the set of all $\mathbb{F}_{q^2}$--rational points of $X$. Let also $P$ be a point in $X(\mathbb{F}_{q^2})$ and let $m_i=m_i(P)$ be a pole number at $P$.

Then according to Lewittes,
\cite[th 1(b)]{Lewittes:90},\cite[p.46 ]{maximal:97}
\[
\#X(\mathbb{F}_{q^2})= N\leq q^2 m_i +1.
\]
Combining the above with the Hasse-Weill bound for a maximal curve, we obtain the following bound
($m_i$ is any pole number)
\[
 \#X(\mathbb{F}_{q^2})=q^2+1 +2g q\leq q^2 m_i +1,
\]
or
\[
 g\leq \frac{q(m_i-1)}{2}.
\]
If $m_i =q$ the above is a result due to  Ihara's, \cite{Ihara} and the
 equality is obtained when $X$ is the Hermitian curve.
Notice that if $P$ is $\mathbb{F}_{q^2}$ rational point then $q,q+1$ are always  pole numbers, \cite[proposition 1.5,(iv)]{maximal:97}.

\section{Hasse-Witt matrix and semigroups}
%
%
%
\label{sec-HasseWitt}
In \cite{Stoehr-Vianna} K.O. St\"ohr and P.Viana introduced a completely local construction of the 
Hasse-Witt matrix \cite{Hasse-Witt}. One of their results that will be useful to us is the following 
\begin{proposition} \label{SV-divi}
Let $(a_{ij})$ be the Hasse-Witt matrix and consider the product
\[
 A_{r}:=(a_{ij})(a_{ij}^p)\cdots (a_{ij}^{p^{r-1}})
\]
{Let $P$ be any point in the curve in question.} 
For each positive $r$ the rank of the matrix $A_r$ is larger than or equal to the 
number of gaps at divisible by $p^r$.   
\end{proposition}
\begin{proof}
 \cite[cor. 2.7]{Stoehr-Vianna}
\end{proof}
Notice that $(a_{ij})$ is dual to the Cartier operator and the matrix $A_{r}$ corresponds
to the application of $r$-times of the corresponding $p$-linear map. For more information 
concerning the Hasse-Witt matrix and the Cartier operator we refer to \cite{BayerGonzalez}.
The rank of  $A_r$  equals  the rank of the Cartier Operator.

Also notice that if the rank is zero, 
then there are no gaps divisible by $p^r$ and every number divisible by $p^r$ 
is a pole number. This can also be seen by different methods, see \cite{RiReJa}.
\begin{remark}
 The $p$-rank of the Jacobian is the rank of the matrix 
$A_g$. Having $p$-rank zero does not give us any information about pole numbers 
since every number greater than $2g$ is a pole number. For this notice that for 
$p\geq 5$  we have  $2g < p^g$. 
\end{remark}

\begin{proposition}\label{ba1}
If a curve has Hasse-Witt matrix zero then every integer divisible by $p$
is a pole number. This implies that $G_1(P)$ is at most an extension of an elementary 
abelian group with a cyclic group of order $p$.
\end{proposition}
\begin{proof}
 If there was a pole number $m$ such that $m<p$ then $G_1(P)$ is 
faithfully represented in $L(mP)$ and it should be elementary abelian. 
If $p<m$, then $m_1=p$ and by lemma \ref{div-inv} we have that
$\ker \rho_1$ divides $p$ therefore $\ker \rho_1$ is either trivial or isomorphic 
to a cyclic group of order $p$. The group $G_1(P)$ is given by a 
short exact sequence
\[
 1 \rightarrow \ker \rho_1 \rightarrow G_1(P) \rightarrow V \rightarrow 1
\]
where $V$ is an elementary abelian group.
\end{proof}

\begin{example} A classical example of a curve with nillpotent Cartier 
operator is given by the Hermitian curve
\[
 y^{p^r}+y=x^{p^r+1},
\]
(which is isomorphic to the Fermat  curve
$x^{p^r+1}+y^{p^r+1}+1=0$), with Cartier operator satisfying  $C^r=0$ \cite{RiReJa}.
\end{example}

\begin{corollary}\label{baoderG_2}
If $X$ has Hasse-Witt matrix zero, and is a big action then 
$G_2(P)$ is cyclic of order $p$.
\end{corollary}
\begin{proof}
If the first pole number is not divisible by $p$, 
then we have a faithful two dimensional representation 
of $G_1(P)$ on $L(m_rP)$, so $G_1(P)$ is elementary abelian, 
therefore bounded by a linear bound on $g$ 
\cite{NakAbel}.
The first pole number is divisible by $p$ since the Hasse-Witt matrix 
is assumed to be  zero. Also the first pole number is  $m_1=|G_2(P)|$ by lemma \ref{big11}. 
Thus, $G_2(P)$ is a cyclic group of order $p$.
\end{proof}

\begin{corollary}\label{nilpotentnonclas}
Let $X$ be a curve of genus $g \geq 2$. 
 If the curve $X$ has nilpotent Cartier operator, i.e., $C^\ell=0$, and moreover 
$p^\ell \leq g$, then 
the curve is non-classical with respect to the canonical  linear series. 
Moreover all  curves with zero Cartier operator that are equipped with an 
automorphism group that has a wild ramified point  and $g\neq p-1$,
 are non-classical with respect to the 
canonical  linear series  with only one hyperelliptic exception, namely $y^2=x^p-x$.
\end{corollary}
\begin{proof}
 Recall that a curve is classical for the canonical divisor if and only if the gap sequence
 for a non Weierstrass point 
is given by $\{1,2,\ldots,g\}$. Proposition \ref{SV-divi} implies that $p^\ell$ is a pole number. 
Since $p^\ell  \leq g$ the 
curve can't be classical. 

Assume now that the curve $X$ has an automorphism group $G$ such that 
there is a wild ramified point in the cover $X \rightarrow X/G$, and has zero 
Cartier operator.  Therefore 
 $\ell=1$, and according to \cite{Roq:70} the existence of wild ramification 
forces  $p-1\leq g$  or the curve is the hyperelliptic curve $y^2=x^p-x$.
For the case $p-1 < g$ the result follows. We don't know what happens for 
the $g=p-1$ case (for small primes, $p<5$, these curves are classical because 
they satisfy the criterion in eq. (\ref{pgrdegK}) below). 

The hyperelliptic curve $y^2=x^p-x$ is a curve of genus  $g=(p-1)/2$
and it is also an Artin-Schreier extension of the rational function field. 
In \cite{ValHyper} it is proved that this curve has zero Cartier operator.
 In fact this is the superspecial hyperelliptic curve with the biggest possible
 genus (see \cite[Theorem 1.2]{Cartierpoints}). 
It is well known that all the hyperelliptic curves of arbritrary characteristic
are $K$--classical, \cite[Satz 8]{Schmidt:39}.
\end{proof}
%

Notice that all superspecial hyperelliptic curves $X$ are K--classical 
and they satisfy the  criterion given in the following equation 
\begin{equation}\label{pgrdegK}
p>\deg K=2g-2,   \Longrightarrow X \textrm{is classical with respect to } K,
\end{equation}
 (see for instance \cite[Theorem 15]{Laksov:84} and \cite[prop. 14.2.64, p.561]{SalvadorBook}),
since from \cite[Theorem 1.2]{Cartierpoints} their genera are upper bounded by $\frac{p-1}{2}$.
Moreover when $\mathrm{char}k=0$ or when eq. (\ref{pgrdegK}) is valid then $X$ is K--classical.

The remarkable fact is that \textit{neither of} the wild ramified coverings
$X\longrightarrow X/G$, with $X$ being \textit{classical or not} with respect to the canonical linear series, 
satisfy  equation \ref{pgrdegK} for $g\geq 3$ (the reader should exclude the hyperelliptic exceptional case considered in corollary \ref{nilpotentnonclas}). This statement
follows from the simple facts that for these curves $p\leq  g+1$, and $2g-2\geq g+1$ for every $g\geq 3$. Keep also in mind that there do not exist non K--classical curves for $g\leq 3$, with only one exception for $p=g=3$, see \cite{Komiya}.

\begin{remark}
Corollary \ref{nilpotentnonclas}, restricted to the the world of maximal curves is  similar 
to the construction 
\cite[proposition 1.7.]{maximal:97}. Indeed,
from theorem 3.3 in \cite{GarciaTafazolian}, every maximal 
and minimal  curve over $\F_{q^2}$, $q=p^\ell$  have nilpotent 
Cartier operator with $C^\ell=0$. The small  difference in the lower bound that is
 given there, $p^\ell-1\leq g$ is explained 
because if $X$ is classical and $g=p^\ell -1$ then $m_1 =p^\ell$ and from 
 \cite[proposition 1.5 and the remark just befor this]{maximal:97} 
this curve should be the 
Hermitian, that  is a contradiction since the Hermitian curve has genus 
$\frac{p^\ell(p^\ell-1)}{2}$.
Observe also that the genus of the nilpotent curves is bounded by
\begin{equation}\label{gbound}
g\leq\frac{p^\ell (p^\ell -1)}{2},
\end{equation}
where $\ell$ is the rank of nilpotency \cite[th. 4.1]{RiReJa}. 
%
\end{remark}

\begin{remark}\label{lastremark}
Combining the results from proposition \cite[Proposition 1.10]{maximal:97}, lemma \ref{llee} here, and  \cite[theorem 2.5]{embedingmaximal},
 we get the following:

Consider a  maximal curve over $\mathbb{F}_{q^2}$, {with  genus $g\geq 2$} and 
the set $\Sigma$ of  $\mathbb{F}_{q^2}$-rational points such 
that the Weierstrass semigroup  up to $q+1$, and hence all the 
Weierstrass semigroup, is generated by two integers.
Then the Weierstrass semigroup at all points of $\Sigma$ is  symmetric,
 i.e. their max gap is always at $2g-1$. 

Since $q,q+1$ must always be in the Weierstrass semigroup at such a
 point, this condition to the numbers of generators,  corresponds to
 the minimum number of generators that a maximal curve can have. Thus we can rephrase:\\
{\em
The maximal curves over $\mathbb{F}_{q^2}$ with minimal set of generators for their
 Weierstrass semigroups at a $\mathbb{F}_{q^2}$--rational point,
have symmetric Weierstrass semigroups at this point.
}
\end{remark}
The next proposition shows that the condition on the number of generators of 
the Weierstrass semigroup in a $\mathbb{F}_{q^2}$-rational point is not 
necessary so that the point has symmetric Weierstrass semigroup. 

\begin{proposition}\label{symGiKu}
Let $X_{GK}$ be the maximal curve over $\mathbb{F}_{q^2}$ defined in
 Giulietti--Korchm{\'a}ros, \cite{newfam}. Then the Weierstrass semigroup at the
 rational point  $X_\infty$ is symmetric. This is an example of curve where the equality in proposition \ref{genera} fails.
\end{proposition}
\begin{proof}

We will use the notation from \cite{newfam}. Let $n=p^h$, $p$ a prime,  $h\geq1$ and 
$q=n^3$. From \cite[equation 10, p. 236]{newfam}, we can write
\begin{eqnarray}\label{maxsymGK}
2g_{GK}-1&=&\sum_{i=1}^3 \biggl(\frac{d_{i-1}}{d_i}-1\biggr)\alpha_i\nonumber\\
&=&-\alpha_1+(n^2-n)\alpha_2+(n-1)\alpha_3,
\end{eqnarray}
where $\alpha_1=n^3-n^2+n,\ \alpha_2=n^3, \ \alpha_3=n^3+1$, are the generators of $H(X_\infty)$, \cite[proposition 5]{newfam}, $d_0=0, \ d_1 = \alpha_1,\ d_2=\mathrm{g.c.d.}(\alpha_1,\alpha_2),$ and $d_3 =\mathrm{g.c.d.}(\alpha_1,\alpha_2, \alpha_3)$.
Suppose now that $2g_{GK}-1$ is a pole number. From \cite[lemma 5]{newfam}
there are uniquely determined no negative integers $j_i$, with $i=1,\ldots 3$ and $j_1,j_2\leq n^2-n$, $j_3\leq n-1$ such that 
\[
2g_{GK}-1=\alpha_1 j_1 +\alpha_2 j_2 +\alpha_3 j_3.
\]
From eq. \ref{maxsymGK} this is equivalent to
\[
(j_1+1)\alpha_1+j_2 \alpha_2+j_3\alpha_3=(n^2-n)\alpha_2+(n-1)\alpha_3
\]
which shows that $j_2=n^2-n, \ j_3=n-1,\ j_1=-1$, a contradiction! 

From \cite[theorem 2]{newfam}, we have that $2g_X=(n^3+1)(n^2-2)+2$. With an immediate calculation
 we show that 
$2g_X <(\alpha_1-1)(\alpha_2-1)$. Indeed, notice that $(\alpha_1,\alpha_2)=n$ so we can not apply lemma 
\ref{llee}. However we can show that $2g_{GK}< (\alpha_1-1)(\alpha_3-1)$ since the right hand of the 
inequality is greater that $(\alpha_1-1)(\alpha_2-1)$. Also $2 g_{GK} \leq (\alpha_2-1)(\alpha_3-1)$ 
since otherwise $X_{GK}$ is Hermitian, a contradiction.
\end{proof}

\begin{remark}
Let $X$ be the Hermitian, or the Garcia--Stichtenoth curve, \cite{maxnotgalois}, or the Giulietti--Korchm{\'a}ros maximal curve over $\mathbb{F}_{q^2}$, $q=p^{h}.$ In fact these are three of the five known families of maximal curves over $\mathbb{F}_{q^2}$ (in the sense that every known maximal curve arise as an $\mathbb{F}_{q^2}$--cover of these curves). The other two families are the Deligne--Lusztig curves that are $\mathbb{F}_{q^2}$ maximal curves with $q$ being a certain power of three and two respectively (see \cite[introduction]{newfam}).  Then $X$  has
 symmetric Weierstrass semigroups in a $\mathbb{F}_{q^2}$ rational point $P$. The Hermitian and the Garcia--Stichtenoth curves  are examples where we obtain the equality in proposition \ref{genera}. 

Indeed, Hermitian and the Garcia--Stichtenoth curves  have  symmetric Weierstrass semigroups in $\mathbb{F}_{q^2}$ rational points because there are
 Artin--Schreier curves,
thus the  Weierstrass semigroup in a $\mathbb{F}_{q^2}$ rational point is generated by 2 elements. From remark \ref{lastremark},  the Weierstrass semigroup
in any rational point is symmetric. For the Hermitian case, we have that the generators
 of $H(P)$ are $d_1=q,\ d_2 =q+1$ and thus $g_H=\frac{q(q-1)}{2}$. 
For the Garcia--Stichtenoth curve $y^{\ell^2}-y=x^{\ell^2-\ell+1}$, $d_1=\ell^2-\ell+1, d_2=\ell^2$ and $g_{GC}=\frac{(d_1-1)(d_2-1)}{2}$, \cite[theorem 1]{maxnotgalois}. The assertion for the Giulietti--Korchm{\'a}ros maximal curve comes from proposition \ref{symGiKu}.  
\end{remark}
\begin{remark}
Let $X$ be  the Hermitian, or the Garcia--Stichtenoth curve,  or the Giulietti--Korchm{\'a}ros 
maximal curve over $\mathbb{F}_{q^2}$. 
Denote by $P\in X$ the $\mathbb{F}_{q^2}$--rational place 
where $H(P)$ is symmetric and with $F$ its function field. 
Take $G$ to be a $p$ subgroup of the automorphism group of $X$. 
Take the Galois cover $\pi : X\longrightarrow X/G$ and suppose 
that there is only one place ramified in that cover. 
We will prove that the  Weierstrass semigroup at $P$ is symmetric 
if and only if the Weierstrass semigroup at the point $\pi(P)$ is symmetric. 

Since every $p$-group is solvable we can decompose the 
cover $\pi$ to a sequence of Artin-Schreier covers.  
So it is enough to consider the case of covers of the form
 $F=F_0(y) $ where
\begin{equation} \label{eqAS1}
 y^{p^h}-y=f, \mbox{ where } f\in F_0,
\end{equation}
and  $f$ has a unique pole at the rational place  $\pi(P)$
of $F_0$, with $v_{P_0}(f)=-m<0$. Then $F/F_0$ is totally ramified  at the place $P$,
 with $P|\pi(P)$.
 Notice that the conditions of Lewittes  theorem, \ref{Lewittes} here, are satisfied. Now we can say the following:
 Lewittes showed that if $g_{X/G}=0$ then $H(P)$ is symmetric;
 If $g_{X/G}>0$
and  $H(\pi(P))$ is symmetric  then $H(P)$ is symmetric, see \cite[p.36 after corollary]{Lewittes:89}.

Now we will prove that if $H(P)$ is symmetric then $H(\pi(P))$ is symmetric. 
According to  Lewittes, theorem \ref{Lewittes} here,
 the max gap in $P$ which equals to $2g_X -1$ is given by
 \begin{eqnarray}\label{Lewsym}
 2g_X -1 &=&m(p^h -1)+p^h h_{g_{X/G}}\nonumber\\
&=&(m+1)(p^h-1)+p^hh_{g_{X/G}}-p^h +1\nonumber\\
&=&\deg\mathrm{Diff}+p^h(h_{g_{X/G}}-1) +1.
 \end{eqnarray}
 Combining eq. (\ref{Lewsym}) and Riemann--Hurwitz formula we should have that $h_{g_{X/G}}=2g_{X/G}-1$, and that means that the Weierstrass semigroup at $\pi(P)$  is also
 symmetric. 

Notice that $\pi(P)$ is also an  $\mathbb{F}_{q^2}$ rational point from the transitivity of the relative degrees.
Finally notice that all the curves considered by Stichtenoth--Garcia--Xing in \cite[Section 3, Theorem 3.2]{GSXcomp} (see also \cite[theorem 2, $m=1$ case]{maxnotgalois}), have symmetric Weierstrass semigroups at $\pi(P)$.
\end{remark}

\def\cprime{$'$}
\providecommand{\bysame}{\leavevmode\hbox to3em{\hrulefill}\thinspace}
\providecommand{\MR}{\relax\ifhmode\unskip\space\fi MR }
\providecommand{\MRhref}[2]{%
  \href{http://www.ams.org/mathscinet-getitem?mr=#1}{#2}
}
\providecommand{\href}[2]{#2}

\end{document}